\documentclass[hidelinks,onefignum,onetabnum]{siamart220329}

\usepackage{graphicx}
\usepackage{epstopdf}
\ifpdf
  \DeclareGraphicsExtensions{.eps,.pdf,.png,.jpg}
\else
  \DeclareGraphicsExtensions{.eps}
\fi

\usepackage{amsmath,amssymb,amsfonts,amscd}
\usepackage{hyperref}
\usepackage[all,cmtip]{xy}
\usepackage{enumerate}
\usepackage{enumitem}
\usepackage{amssymb, latexsym, amsmath}
\usepackage{url}
\usepackage{nicematrix}
\usepackage{tikz}
\usepackage{algpseudocode,algorithm}
\usepackage{subcaption}
\usepackage{mathtools}
\usepackage[normalem]{ulem}

\usepackage[normalem]{ulem}

\newcommand{\GN}{generalized Nyström}
\newcommand{\GNa}{generalized Nyström approximation}
\newcommand{\HMTa}{HMT approximation}
\newcommand{\RR}{Rayleigh-Ritz}
\newcommand{\tU}{\Tilde{U}}
\newcommand{\tV}{\Tilde{V}}
\newcommand{\bA}{\Bar{A}}
\newcommand{\bbA}{\Bar{\Bar{A}}}
\newcommand{\tA}{\Tilde{A}}
\newcommand{\Hp}{H_{p}}
\newcommand{\Hphat}{\hat{H}_{p}}
\newcommand{\Fp}{F_{p}}

\newcommand{\btau}{\Bar{\tau}}

\newcommand{\minindex}{k}

\newcommand{\Nr}{N_r}
\newcommand{\tNr}{\tilde{N}_r}
\newcommand{\calO}{\mathcal{O}}

\newcommand\mtiny[1]{\mbox{\tiny\ensuremath{#1}}}
\newcommand*{\Scale}[2][4]{\scalebox{#1}{\ensuremath{#2}}}%

\newsiamremark{remark}{Remark}
\newsiamremark{hypothesis}{Hypothesis}
\crefname{hypothesis}{Hypothesis}{Hypotheses}
\newsiamthm{claim}{Claim}

\headers{{Extracting accurate singular values from subspaces}}{L. Lazzarino, H. Al Daas, Y. Nakatsukasa}

\title{Matrix Perturbation Analysis of Methods for Extracting Singular Values from Approximate Singular Subspaces}

\author{Lorenzo Lazzarino\thanks{Mathematical Institute, University of Oxford, Oxford, OX2 6GG, UK  (\email{lorenzo.lazzarino@maths.ox.ac.uk, nakatsukasa@maths.ox.ac.uk}).}, 
\and Hussam Al Daas\thanks{Scientific Computing Department, STFC, Rutherford Appleton Laboratory, Harwell Campus, Didcot, Oxfordshire, OX11 0QX, UK  (\email{hussam.al-daas@stfc.ac.uk}).}
 \and Yuji Nakatsukasa\footnotemark[1]}

\usepackage{amsopn}

\begin{document}

\maketitle

\begin{abstract}
    Given (orthonormal) approximations $\tU$ and $\tV$ to the left and right subspaces spanned by the leading singular vectors of a matrix $A$, we discuss methods to approximate the leading singular values of $A$ and study their accuracy. In particular, we focus our analysis on the \GNa{}, as surprisingly, it is able to obtain significantly better accuracy than classical methods, namely \RR{} and (one-sided) projected SVD. 
    A key idea of the analysis is to view the methods as finding the exact singular values of a perturbation of $A$. In this context, we derive a matrix perturbation result that exploits the structure of such $2\times2$ block matrix perturbation. This leads to sharp bounds that predict well the approximation error trends and explain the difference in the behavior of these methods. Finally, we present an approach to derive an a-posteriori version of those bounds, which are more amenable to computation in practice.
\end{abstract}
% REQUIRED

% REQUIRED
%%%%%%%%%%%%%%%%%%%%%%%%%%%%%%%%%%%%%%%%%%%%%%%%%%%%%%%%%%%%%%%%%%%%%%%%%%%%%%%%%%%%%%%%%%%%%%%%%%%%%%%%%%%%%%%%%%%%%%%%%%%%%%%%%%%%%%%%%%%%%%%%
%%%%%%%%%%%%%%%%%%%%%%%%%%%%%%%%%%%%%%%%%%%%%%%%%%%%%%%%%%%%%%%%%%%%%%%% INTRODUCTION %%%%%%%%%%%%%%%%%%%%%%%%%%%%%%%%%%%%%%%%%%%%%%%%%%%%%%%%%%
%%%%%%%%%%%%%%%%%%%%%%%%%%%%%%%%%%%%%%%%%%%%%%%%%%%%%%%%%%%%%%%%%%%%%%%%%%%%%%%%%%%%%%%%%%%%%%%%%%%%%%%%%%%%%%%%%%%%%%%%%%%%%%%%%%%%%%%%%%%%%%%%
\section{Introduction} \label{sec:Intro}
Approximating the singular values of a matrix $A$ is important for several reasons across various fields of mathematics, engineering and data science. In principal component analysis (PCA), the singular values of the data matrix provide information on the variance captured by each principal component. In linear inverse problems, access to singular values is important to build effective regularisation techniques \cite{Han10,HanPS13,Ste93}. Additionally, estimating leading singular values is essential for determining where to truncate to achieve a meaningful low-rank approximation. This is beneficial in any context where low-rank approximations are applied.

Subspace iteration, Golub-Kahan bi-diagonalization, and randomized techniques are widely used methods to approximate the left and right dominant singular subspaces. Once an approximation to the left/right singular subspaces $\tU$ and/or $\tV$ is available, a number of strategies exist to extract the singular values. The Rayleigh-Ritz (RR) technique \cite{daxExtremumPropertiesOrthogonal2010,saadNumericalMethodsLarge2011,xin-guoRayleignQuotientTheory1992} is the standard method used to extract the singular values once the singular subspaces are approximated. 
A significant motivation behind this work is the efficient randomized methods for low-rank approximation, such as the Halko-Martinsson-Tropp (HMT) \cite{halkoFindingStructureRandomness2011} and the generalized Nystr{\"o}m (GN) \cite{clarksonNumericalLinearAlgebra2009,nakatsukasaFastStableRandomized2020,woolfeFastRandomizedAlgorithm2008}, which are based on first finding approximations to the left and/or right dominant singular subspaces of $A$.

Motivated by streaming data analysis, RR, HMT, GN, and other methods to extract an approximation of the singular values can be categorized into two families, one-pass and multi-pass, based on the number of accesses to the matrix of interest that is required to obtain the approximation.
It can be observed numerically that these methods have different approximation behavior. 
The goal of this paper is to understand and characterize the reasons for such difference. This is accomplished by interpreting the methods as perturbations of the original matrix and studying their accuracy through a certain matrix perturbation result, to obtain error bounds for individual singular values.

By design, GN and HMT are (randomized) algorithms for finding a low-rank approximation $A\approx \tilde A$. 
The approximation quality of GN has been extensively analyzed in \cite{tropp2017practical,nakatsukasaFastStableRandomized2020}, showing roughly that the error is optimal up to a modest factor and a slight oversampling. Analogous results are available for HMT, which comes with a slightly better error~\cite[\S~10]{halkoFindingStructureRandomness2011}. These results hold with very high quantifiable probability. It is possible to use such matrix approximation bounds to derive error bounds for the singular values. For example, a low-rank approximation error bound in the 2-norm 
$\|A-\tilde A\|_2$
would give an error bound for all singular values. However, as we will see, the errors are often not uniform, with the leading singular values typically having significantly higher accuracy. Such property cannot be captured by a global, low-rank approximation error. 
For HMT (and more generally subspace iteration), convergence analysis of singular values is provided in~\cite{gu2015subspace,saibaba2019randomized}.

We review in \Cref{sec:background} several methods to approximate the leading singular values of a matrix given the approximate singular subspaces. Motivated by numerical observation showing its better accuracy compared to other single-pass strategies, we analyze the \GNa{} and generalize the analysis to the \RR{} and (randomized) SVD approximations. Specifically, in \Cref{sec:GNPerturbation}, we present a formal interpretation of the \GNa{} as a structured matrix perturbation and detail its underlying structure. Classic results from perturbation analysis of singular values such as Weyl’s inequality provide pessimistic bounds that cannot provide insights into the expected accuracy of the approximated singular values by each of these methods. Therefore, in \Cref{sec:MatPert} we develop results that exploit the structure of the perturbation. In particular, we study perturbations of  $2\times 2$ block matrices. This result is a generalization of the result for eigenvalues based on studying perturbations of symmetric matrices in \cite{nakatsukasaEigenvaluePerturbationBounds2012}. Then, in \Cref{sec:ApplToMethods} we apply this formally derived and general result to the problem above, obtaining bounds for each of the mentioned strategies. Finally, in \Cref{sec:NumIll}, we illustrate how these bounds are often sharper than Weyl's inequality, how they can be used to compare methods, and we present an a-posteriori alternative to them.

From a computational perspective, our results suggest the following strategy for approximating singular values. 
When only one subspace $\tilde U$ or $\tilde V$ is available, the SVD approach is recommended. When $\tilde U,\tilde V$ are both available, \GN{} tends to give significantly higher accuracy than other single-pass methods, as our bounds indicate and Figure~\ref{fig:Motivation} illustrates; we regard this as a key contribution of this paper. When multiple passes with $A$ are allowed, we recommend HMT, or better yet, combined with several iterations of subspace iteration as much as the computational budget permits.

\section{Background and Motivation}\label{sec:background}
Let $A \in \mathbb{R}^{m\times n}$ and its economic singular value decomposition $A = U\Sigma V^*$ with singular values $\sigma_1 \ge \cdots \ge \sigma_n \ge 0$.
The first singular values and their corresponding singular vectors are often referred to as \textit{leading}, and they play a crucial role in applications involving e.g. low-rank approximation. 
%In practice, these are often the most important part of the decomposition. 
Throughout the paper, we say that a matrix $M$ is orthonormal if it is a tall matrix with orthonormal columns, i.e., $M^*M = I$. For simplicity, we assume the matrices are real, although everything carries over verbatim to the complex case. Thus, the superscript $*$ denoting the conjugate transpose can be interpreted as the matrix transpose.

Suppose we have access to orthonormal approximations $\tU$ and $\tV$ to the leading singular subspaces of the matrix $A$. Sometimes we have approximations to different dimensions of left and right singular subspaces, that is, $\tV \in \mathbb{R}^{n\times r}$ and $\tU \in \mathbb{R}^{m\times (r + \ell)}$, orthonormal approximations to the leading $r$-dimensional left and $(r + \ell)$-dimensional right singular subspaces, where $\ell \geq 0$ is the so-called \textit{oversampling parameter}. In other cases, we may have just one of them. In either case, a classical numerical linear algebra problem is to approximately find the leading $r$ singular values of $A$. There are multiple strategies to accomplish such a task. The most classical one is the \textit{\RR{} method}, a strategy usually presented for the extraction of eigenvalues \cite[Sec. 4-5]{saadNumericalMethodsLarge2011}, but that can be easily generalized to singular values \cite{daxExtremumPropertiesOrthogonal2010,xin-guoRayleignQuotientTheory1992}. It extracts the leading $r$ singular values by considering the singular values of the matrix $A_{\mtiny{RR,\tV,\tU}}:=\tU^*A\tV$, i.e., for $i=1, \dots,r$,
\begin{displaymath}
  \sigma_i \approx \sigma_i(\tU^*A\tV) =:\sigma_i^{\mtiny{RR}}.
\end{displaymath}
If we have access to only one of the two approximations, say $\tV$, then we can extract the leading singular values by looking at the matrix $A_{\mtiny{ SVD,\tV}}:= A\tV$, i.e., for $i=1, \dots,r$,
    \begin{displaymath}
      \sigma_i \approx \sigma_i(A\tV) =:\sigma_i^{\mtiny{SVD}}.
    \end{displaymath}
We refer to this as the \textit{SVD approximation}, as it is based on a (one-sided) projected SVD. 
%\rr{For example, the singular values computed by HMT (below) can be seen as the output of SVD from the left singular subspace $\tilde U$, $\sigma_i(\tU^*A) $.}
Less intuitive methods can also be employed. The general idea is to use the available approximations in the expression of low-rank approximations of $A$. For example, we can consider the \GN{} low-rank approximation \cite{clarksonNumericalLinearAlgebra2009,nakatsukasaFastStableRandomized2020,woolfeFastRandomizedAlgorithm2008}. That is, for some $X\in \mathbb{R}^{n\times r}$ and $Y\in \mathbb{R}^{m\times (r+\ell)}$, we construct an approximation of $A$ with column space $AX$ and row space $Y^*A$. In particular,
\begin{equation}
\label{def:GN}
    A_{\mtiny{ GN,X,Y}}:= AX(Y^*AX)^\dagger Y^*A,
\end{equation}
where the superscript $\dagger$ denotes the Moore–Penrose pseudoinverse as defined in \cite[Sec. 5.5.2]{golubMatrixComputations2013}. This expression has been first introduced in \cite{wedderburnLecturesMatrices1934} and then used and analyzed in the randomized low-rank approximation context in \cite{clarksonNumericalLinearAlgebra2009,nakatsukasaFastStableRandomized2020,woolfeFastRandomizedAlgorithm2008}. In its classic definition, $X$ and $Y$ are taken to be random matrices, e.g., Gaussian matrices or more structured matrices such as SRTT \cite{boutsidisImprovedMatrixAlgorithms2013,troppIMPROVEDANALYSISSUBSAMPLED2011}. Thus, we can interpret this approximation as a two-sided projection process. The idea is to use the \GN{} expression to obtain an approximation to singular values. That is, we use $\tU$ and $\tV$ in \cref{def:GN} (instead of $X$ and $Y$) and then compute the singular values, i.e., for $i=1, \dots,r$,
\begin{displaymath}
  \sigma_i \approx \sigma_i(A\tV(\tU^*A\tV)^\dagger\tU^*A) =:\sigma_i^{\mtiny{GN}}.
\end{displaymath}
We will refer to this as the \textit{\GNa}, see \Cref{algo:GN} for an implementation. It is worth noting that, unlike the other single-pass methods, $\tU$ and $\tV$ need not be orthonormal in \GN. We will however assume they are, as this simplifies the analysis without losing generality.

A similar approach can be used for the case where just $\tV$ is available. We consider the low-rank approximation computed by the \textit{Halko-Martinsson-Tropp (HMT) method}, also known as randomized SVD \cite{clarksonLowRankApproximationRegression2017,halkoFindingStructureRandomness2011,rokhlin2009randomized}. 
Given a random matrix $X\in \mathbb{R}^{n\times r}$, we compute the product $AX$ and perform a QR factorization $AX = QR$. Then, we define the low-rank approximation by
\begin{equation}
    \label{def:HMT}
    A_{\mtiny{ HMT,X}} = Q(Q^*A).
\end{equation}
The idea is to use the HMT process to extract singular values. Thus, we start the HMT process with $\tV$ (instead of $X$) and then compute the singular values of $Q^*A$. We have, for $i=1, \dots,r$,
\begin{displaymath}
  \sigma_i \approx \sigma_i(Q^*A) =:\sigma_i^{\mtiny{HMT}}.
\end{displaymath}
Note, that due to the orthogonality of $Q$, there is no need to perform the left multiplication by it in \cref{def:HMT} for the purpose of extracting singular values. We will refer to this as the \textit{\HMTa{}}, see \Cref{algo:HMT} for an implementation.

\begin{figure}[t]
    \begin{subfigure}{.5\textwidth}
          \centering
          \includegraphics[width=\linewidth]{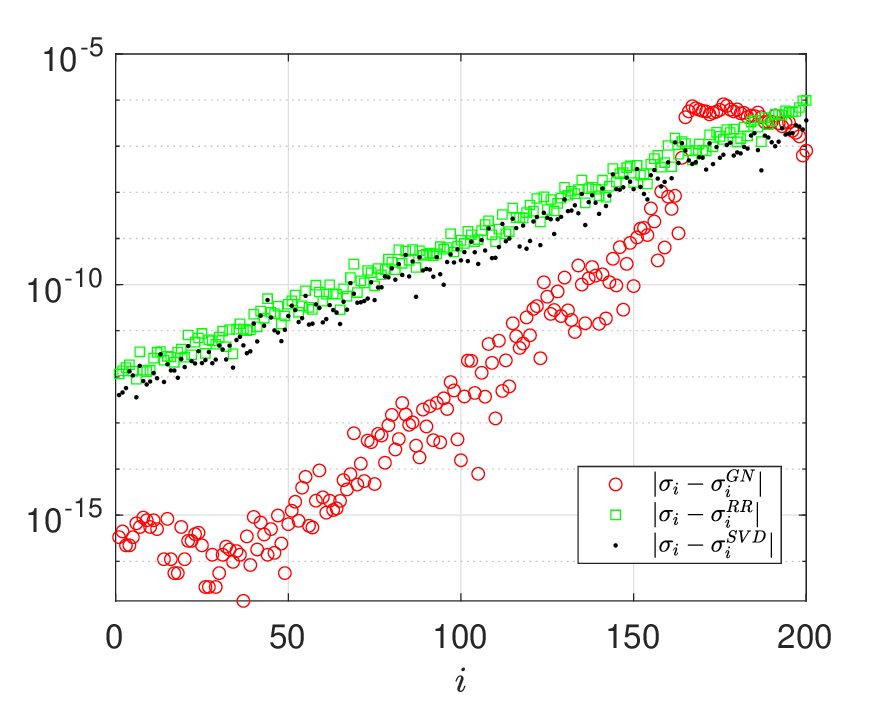}
          \caption{Without Oversampling ($\ell = 0$)}
          \label{fig:Motivation_exp}
    \end{subfigure}%
    \begin{subfigure}{.5\textwidth}
          \centering
          \includegraphics[width=\linewidth]{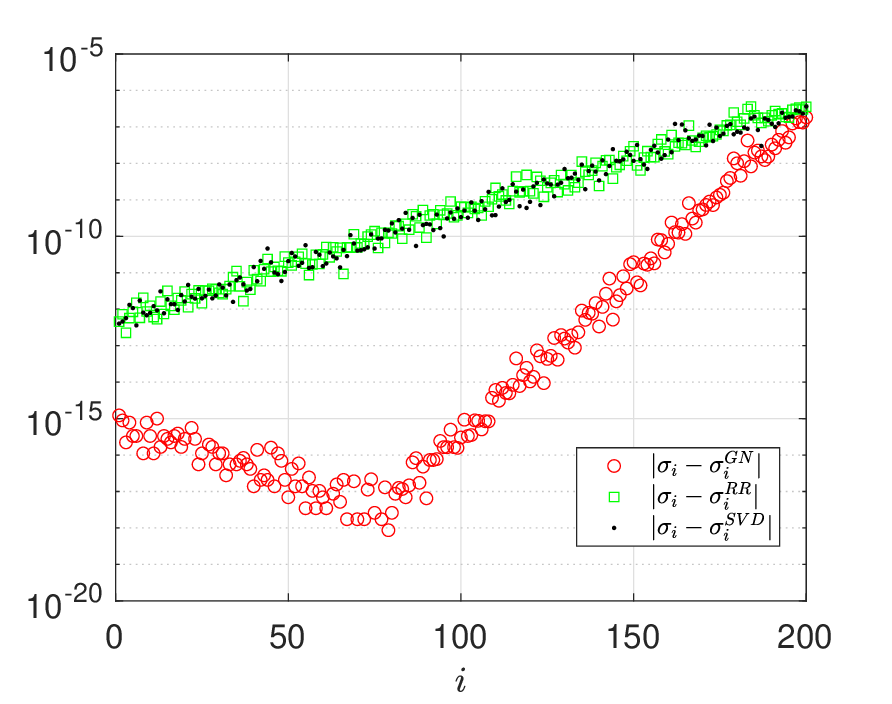}
          \caption{With Oversampling ($\ell > 0$)}
          \label{fig:Motivation_exp_ov}
    \end{subfigure}
    \caption{\textbf{Accuracy of single-pass extraction methods.} We compare the approximation quality of the studied single-pass methods: \GN{} (red circles), \RR{} (green squares), and SVD (black dots). The matrix in this example is constructed such that its singular values decay exponentially. 
The V-shaped errors with \GN\ is likely an artifact of roundoff errors; in exact arithmetic, we expect the error to be roughly proportional to $1/\sigma_i$. Further details on the experiment are given in \Cref{sec:NumIll}.}
    \label{fig:Motivation}
\end{figure}
 Although both the HMT and the \GN{} approximations need two multiplications by $A$, only the \GNa{} is single-pass, that is, can be computed by accessing $A$ once (as are the SVD approximation and the \RR{} method). Moreover, all these methods return exact singular values if $\tU$ and $\tV$ are exact leading singular subspaces.
 
 \Cref{fig:Motivation} shows how the \GNa{} typically tends to give the best approximation among single-pass methods, providing significantly better approximations of the singular values, especially for the leading ones. This serves as a motivation to investigate further how to characterize the accuracies of the various strategies to understand why the \GNa{} performance is often significantly better. For this reason, we focus our analysis on the \GNa, and then, in \Cref{subsec:Appltoother}, generalize it to the other methods. 

All these strategies can be interpreted as perturbations of the original matrix $A$ \cite{troppRandomizedAlgorithmsLowrank2023} that, after particular orthogonal transformations, have peculiar structures that we will exploit. A standard result in this field is \textit{Weyl's inequality} \cite[Cor. 7.3.5]{hornMatrixAnalysis2018},\cite[Cor. I.4.31]{stewartMatrixAlgorithms1998}: for any matrix $M$ and for all $i$,
\begin{equation}\label{eq:Weyl}
    |\sigma_i(M) - \sigma_i(M + E)| \leq \|E\|_2.
\end{equation}
An analogous result holds for eigenvalues of symmetric matrices \cite[Cor. 6.3.4]{hornMatrixAnalysis2018},\cite[Thm. 10.3.1]{parlettSymmetricEigenvalueProblem1998},\cite[Cor. III.4.10]{stewartMatrixPerturbationTheory1990}. Both results are usually referred to as Weyl's inequality, bound, or theorem. As is clear from \cref{eq:Weyl}, this bound neither captures nor exploits the structure of the perturbation $E$. Moreover, it does not depend on $i$ so it does not predict the difference between the change of the leading and smallest eigen/singular values. We note that for some structured cases, better bounds than \eqref{eq:Weyl} are available, i.e. those that give relative perturbation bounds \cite{drmacNumericalMethodsAccurate2021,Li1998relativeperturbation}. However, typically relative perturbation bounds are helpful more for the smallest singular values, and not for the leading ones. Notable exceptions are the bounds for the extraction of (leading) singular values by the HMT method in \cite{gu2015subspace,saibaba2019randomized}. For this reason, in \Cref{subsec:Saibaba}, we compare our results for HMT with the one in \cite{saibaba2019randomized}.
%%%%%%%%%%%%%%%%%%%%%%%%%%%%%%%%%%%%%%%%%%%%%%%%%%%%%%%%%%%%%%%%%%%%%%%% GN AS A PERTURBATION OF THE ORIGINAL MATRIX %%%%%%%%%%%%%%%%%%%%%%%%%%%%%%%%%%%%%%%%%%%%%%%%%%%%%%%%%%%%%%%%%%%%%%%%
\section{Generalized Nyström and Matrix Perturbation} \label{sec:GNPerturbation}
The main results of this paper are based on a matrix perturbation interpretation of the aforementioned methods. The following illustrates how the \GNa{} can be interpreted as a perturbation of the original matrix $A$ \cite{troppRandomizedAlgorithmsLowrank2023}. 

For any orthonormal $\tU$ and $\tV$, considering their orthogonal complements $\tV_{\perp}$ and $\tU_{\perp}$, we can define the orthogonal matrices
\begin{displaymath} 
Q_{1}= \begin{bmatrix} \tU & \tU_{\perp}\end{bmatrix}\in\mathbb{R}^{m\times m}, \quad Q_{2}= \begin{bmatrix} \tV & \tV_{\perp}\end{bmatrix}\in\mathbb{R}^{n\times n},
 \end{displaymath}
for which 
 \begin{displaymath}
 \begin{aligned}
     (Q_1^*AQ_2)_{\mtiny{ GN,\begin{bmatrix}I_r \\0\end{bmatrix},\begin{bmatrix}I_{r + \ell} \\ 0\end{bmatrix}}} &=Q_1^*AQ_2\begin{bmatrix}
        I_{r} \\ 0
    \end{bmatrix} \left(\begin{bmatrix}
        I_{r + \ell} & 0
    \end{bmatrix}Q_1^*AQ_2\begin{bmatrix}
        I_{r} \\ 0
    \end{bmatrix}\right)^\dagger\begin{bmatrix}
        I_{r + \ell} & 0
    \end{bmatrix}Q_1^*AQ_2 \\[1ex]
    &=Q_1^* \left( A\tV(\tU^*A\tV)^\dagger\tU^*A\right) Q_2 =Q_1^* \left( A_{\mtiny{ GN,\tV,\tU}}\right) Q_2.
 \end{aligned}
 \end{displaymath}
Due to the orthogonality of the transformations $Q_1$ and $Q_2$, this trivially implies that
 \begin{equation}
 \label{eq:transfSingVal_eq_SingVal}
 \begin{aligned}
 | \sigma_i(A) - \sigma_i(A_{\mtiny{ GN,\tV,\tU}}) | &= | \sigma_i(Q_1^*AQ_2) - \sigma_i(Q_1^*A_{\mtiny{ GN,\tV,\tU}}Q_2) | \\
 &= | \sigma_i(Q_1^*AQ_2) - \sigma_i(\left( Q_{1}^{*}AQ_{2} \right)_{\mtiny{ GN,\begin{bmatrix}I_r \\0\end{bmatrix},\begin{bmatrix}I_{r + \ell} \\ 0\end{bmatrix}}}) |.     
 \end{aligned}
 \end{equation}
This means that for any orthonormal $\tU$ and $\tV$ we can find a set of coordinates in which we can study the \GNa{} with sampling matrices
\begin{displaymath}
\label{ass:identities}
\begin{bNiceMatrix}[first-row, first-col, cell-space-top-limit=3pt, extra-margin=5pt] 
        &r              \\
    r   &I_r        \\
    n-r   &0
\end{bNiceMatrix}, \quad 
\begin{bNiceMatrix}[first-row, first-col, cell-space-top-limit=3pt, extra-margin=5pt] 
        &r+\ell              \\
    r+\ell   &I_{r+\ell}        \\
    m-(r+\ell)   &0 
\end{bNiceMatrix}, 
\end{displaymath} 
applied to the matrix $\bA := Q_1^*AQ_2$. Dividing $\bA$ into blocks of appropriate sizes, i.e., 
\begin{equation}
\label{def:Abar}
\bA :=
\begin{bNiceMatrix}[first-row, first-col, cell-space-top-limit=3pt, extra-margin=5pt] % code-for-first-row = $\boldmath$, code-for-first-col = $\boldmath$
        &r              &n-r                          \\
    r+\ell   &\bA _{11}         &\bA _{12}         \\ 
    m-(r+\ell)  &\bA_{21}    &\bA_{22}            
    
\end{bNiceMatrix},
\end{equation} 
we can analyze the effect the \GNa{} has on each block of $\bA$. We have
\begin{displaymath}
\label{eq:GNgeneralBlocks}
\bA_{\mtiny{ GN,\begin{bmatrix}I_r \\0\end{bmatrix},\begin{bmatrix}I_{r + \ell} \\ 0\end{bmatrix}}} = 
\begin{bmatrix}
    \bA_{11}       \\[1ex]
 \bA_{21} 
\end{bmatrix} \, 
(\bA_{11})^\dagger \, 
\begin{bmatrix}
 \bA_{11}  &   \bA_{12}
\end{bmatrix}.
\end{displaymath}
That is,
\begin{displaymath}
    \label{eq:GNIdentities}
  \bA_{\mtiny{ GN,\begin{bmatrix}I_r \\0\end{bmatrix},\begin{bmatrix}I_{r + \ell} \\ 0\end{bmatrix}}}  = \begin{bmatrix}
    \bA_{11}\bA_{11}^\dag \bA_{11}  & \bA_{11}\bA_{11}^\dag \bA_{12}\\[1ex]
    \bA_{21}\bA_{11}^\dag \bA_{11}  & \bA_{21}\bA_{11}^\dag \bA_{12} 
\end{bmatrix}. 
\end{displaymath} 
A reasonable and generic assumption to make is that the sub-matrix $\bA_{11}$ is full-rank, i.e., $\mbox{rank}(\bA_{11})=r$. In this case, $\bA_{11}$ has full column rank and, thus, the pseudoinverse acts as a left inverse, i.e., $\bA_{11}^\dag \bA_{11} = I_{r}$, \cite[Thm. III.1.2]{stewartMatrixPerturbationTheory1990}. For general $\ell$, the same cannot be said for $\bA_{11} \bA_{11}^\dag$. This implies that, under the above assumption,
\begin{equation}
\label{eq:GNIdentities-final}
    \bA_{\mtiny{ GN,\begin{bmatrix}I_r \\0\end{bmatrix},\begin{bmatrix}I_{r + \ell} \\ 0\end{bmatrix}}}
= \begin{bmatrix}
    \bA_{11} & \bA_{11}\bA_{11}^\dag \bA_{12}\\[1ex]
    \bA_{21} & \bA_{21}\bA_{11}^\dag \bA_{12} 
\end{bmatrix}.
\end{equation}
\Cref{eq:GNIdentities-final} leads to a clear relation between the matrix $A$ and $A_{\mtiny{GN,\tV,\tU}}$, i.e.,
\begin{equation}
\label{eq:GNperturbation}
Q_{1}^{*}\left(A - A_{\mtiny{ GN,\tV, \tU}}\right)Q_{2} = 
\begin{bmatrix}
    0 & \bA_{12} - \bA_{11}\bA_{11}^\dag \bA_{12}\\[1ex]
    0 & \bA_{22} - \bA_{21}\bA_{11}^\dag \bA_{12}
\end{bmatrix} =:E_{\mtiny{GN}}.
\end{equation}
Therefore, the \GNa{} of a matrix $A$ can be seen as a perturbation of $A$ that, up to orthogonal transformations, has the particular structure as in \cref{eq:GNperturbation}, where the non-zero perturbations are, respectively, the difference between $\bA_{12}$ and its orthogonal projection into the column space of $\bA_{11}$ \cite[Thm. III.1.3]{stewartMatrixPerturbationTheory1990} ($(1,2)$ block in \cref{eq:GNperturbation}), and the (generalized) Schur complement \cite{zhangSchurComplementIts2005} of $\bA$ ($(2,2)$ block in \cref{eq:GNperturbation}). In the case where $\tU$ and $\tV$ are given with the same column size (i.e., without oversampling, $\ell =0$), the pseudoinverse acts also as a right inverse, as, in this case, $A_{11}$ has full row rank, and, therefore, the top-right block is also $0$. At first sight, this seems an advantageous situation, however, it is important to notice that the considered blocks change size depending on the size of $\tU$ and $\tV$, resulting in a different, and potentially larger, bottom-right perturbation. 

It is also worth noting that by transforming $A$ into $\bA$ we typically (assuming $\tU, \tV$ are high quality) construct a matrix with small off-diagonal blocks and leading singular values far from the spectrum of the $\bA_{22}$ block. More precisely, the magnitude of the off-diagonal blocks will depend on the accuracy of the approximate subspaces. In \Cref{sec:NumIll}, we see that, for the \GNa, when $A$ is low rank, little information about $A$ can be sufficient to give good enough approximate singular subspaces. In the next section, we aim to derive a bound that exploits both the structure of $\bA$ and the perturbation in \cref{eq:GNperturbation}.

%%%%%%%%%%%%%%%%%%%%%%%%%%%%%%%%%%%%%%%%%%%%%%%%%%%%%%%%%%%%%%%%%%%%%%%%%%%%%%%%%%%%%%%%%%%%%%%%%%%%%%%%%%%%%%%%%%%%%%%%%%%%%%%%
%%%%%%%%%%%%%%%%%%%%%%%%%%%%%%%%%%%%%%%%%%%%%%%%%%%%%%%%% MATRIX PERTURBATION RESULTS %%%%%%%%%%%%%%%%%%%%%%%%%%%%%%%%%%%%%%%%%%
%%%%%%%%%%%%%%%%%%%%%%%%%%%%%%%%%%%%%%%%%%%%%%%%%%%%%%%%%%%%%%%%%%%%%%%%%%%%%%%%%%%%%%%%%%%%%%%%%%%%%%%%%%%%%%%%%%%%%%%%%%%%%%%%
\section{Matrix Perturbation Results} \label{sec:MatPert}
In this section, we step away from the problem of extracting singular values and derive a general matrix perturbation result. We derive a bound on the change of singular values when a perturbation is applied that captures and exploits the structure of the perturbation and that predicts the difference in change between different singular values. Our analysis starts from a result for eigenvalues of symmetric block tridiagonal matrices presented in \cite{nakatsukasaEigenvaluePerturbationBounds2012}. We start by considering the result proven in \cite{nakatsukasaEigenvaluePerturbationBounds2012} and, using a technique based on the Jordan-Wielandt Theorem presented in \cite{liNoteEigenvaluesPerturbed2005}, we derive the result for the general case. 
%%%%%%%%%%%%%%%%%%%%%%%%%%%%%%%%%%%%%%%%%%% 2X2 BLOCK MATRIX %%%%%%%%%%%%%%%%%%%%%%%%%%%%%%%%%%%%%%%%%%%%%%%%%%%%%%%%%%%%%%%%%%

Theorem 3.2 in \cite{nakatsukasaEigenvaluePerturbationBounds2012} gives us a bound on the change in the eigenvalues when a symmetric perturbation is applied to a symmetric $2\times 2$ block matrix. In more detail, considering the symmetric $2\times 2$ block $n\times n$ matrix 
\begin{displaymath}
A := \begin{bmatrix}
    A_{11} & A_{21}^* \\ A_{21} & A_{22}
\end{bmatrix},\end{displaymath}
with eigenvalues $\{\lambda_i\}_{i=1}^{n}$ and the perturbed matrix  
\begin{displaymath}
\hat{A} := A + \begin{bmatrix}
    E_{11} & E_{21}^* \\ E_{21} & E_{22}
\end{bmatrix} =: A + E,
\end{displaymath}
with eigenvalues $\{\hat{\lambda}_i\}_{i=1}^{n}$, we can define, for $i=1,\dots, n$,
\begin{displaymath}
    \tau_i = \left( \frac{\|A_{21}\|_2+\|E_{21}\|_2}{\displaystyle \min_{\minindex}|\lambda_i - \lambda_\minindex(A_{22})| - 2\|E\|_2}\right).
\end{displaymath}
Then, for all $i$ for which $\tau_i > 0$, we have
\begin{equation}
\label{eq:Yujion2x2}
    | \lambda_i - \hat{\lambda}_i | \leq \|E_{11}\|_2 + 2\|E_{21}\|_2 \tau_i + \|E_{22} \|_2 \tau_i^2.
\end{equation}
This bound is derived by considering the eigendecompositions of $A+tE$ for $t\in[0,1]$ and observing that small components of the $t$-dependent eigenvector in positions corresponding to the dominant elements of $E$ implies the eigenvalue has a small derivative with respect to $t$.
Hence, bound \cref{eq:Yujion2x2} is found by bounding those components ( a special case of the Davis--Kahan generalized $\sin\Theta$ theorem \cite[Thm. 6.1]{davisRotationEigenvectorsPerturbation1970} is used) and integrating. The condition $\tau_i >0$ comes from the well-definiteness of the intermediate bound on the components of the $t$-dependent eigenvector and it means that the best situation for this bound is when the analyzed eigenvector is well-separated from the spectrum of $A_{22}$ and the perturbation is small in norm. For \cref{eq:Yujion2x2} to be better than Weyl's inequality for eigenvalues, it is necessary that $\tau_i < 1$.  If $\|E_{11}\| \ll \|E\|$ and $\lambda_i$  is far from the spectrum of $A_{22}$  then $\tau_i \ll 1$. 

In the following, we use a technique presented in \cite{liNoteEigenvaluesPerturbed2005} to show that an analogous result holds for the singular values of a non-symmetric rectangular $2\times 2$ block $m\times n$ matrix 
\begin{displaymath}
    H := \begin{bmatrix}
        G_1 & B \\ C & G_2
    \end{bmatrix}.
\end{displaymath}
This is based on the \textit{Jordan-Wielandt theorem}, \cite[Thm. I.4.2]{stewartMatrixPerturbationTheory1990}, \cite[Thm. 7.3.3]{hornMatrixAnalysis2018}: Let $\{\sigma_i\}_{i=1}^n$ be the singular values of a matrix $M \in \mathbb{R}^{m\times n}$, with $m \geq n$. Then, the symmetric matrix 
\begin{displaymath}
    \begin{bmatrix}
        0 & M \\ M^* & 0
    \end{bmatrix} \in \mathbb{R}^{(m+n)\times (m+n)}
\end{displaymath}
has eigenvalues (in descending order)
\begin{equation}
\label{eq:eigenJW}
\{ \lambda_j \}_{j=1}^{m+n} = \{ \sigma_{1}, \dots , \sigma_{n}, \underbrace{0, \dots, 0}_{m-n \, \text{times}},-\sigma_{n},\dots, -\sigma_1 \}.
\end{equation}
This can be used to study the eigenvalues of 
\begin{displaymath}
   H_{JW} := \begin{bmatrix}
       &0 & | & H &\\ &- & - & - &\\ &H^* & | & 0&
   \end{bmatrix} = \begin{bmatrix}
0 & 0 & | & G_1 & B \\
0 & 0  & | & C & G_2 \\
- & - & - & - & - \\
G_1^* & C^* & | &  0 & 0 \\
B^* & G_2^* & | &  0 & 0 \\
\end{bmatrix}
\end{displaymath}
to obtain results on the singular values of $H$.  However, we still want to obtain results that are translatable back to the original blocks $G_1, G_2, B$, and $C$. Hence, we consider a permutation of $H_{JW}$ with the same eigenvalues but different structure. In other words, we first exchange the second and third block rows (red) and then the second and third block columns (blue) of $H_{JW}$, to obtain the permuted matrix $\Hp$:
\begin{equation}
\Scale[0.8]{
\label{eq:2by2permutations}
\begin{bNiceMatrix}
0 & 0 & | & G_1 & B \\
\Block[fill=red!15,rounded-corners]{1-5}{}
0 & 0  & | & C & G_2 \\ 
- & -  & - & - & - \\
\Block[fill=red!15,rounded-corners]{1-5}{}
G_1^* & C^* & | &  0 & 0 \\
B^* & G_2^* & | &  0 & 0 
\end{bNiceMatrix}  \stackrel{P_{\text{red}}}{\looparrowright} 
\begin{bNiceMatrix}
0 &\Block[fill=blue!15,rounded-corners]{5-1}{} 0 & | &\Block[fill=blue!15,rounded-corners]{5-1}{} G_1 & B \\
G_1^* & C^* & | & 0 & 0 \\
- & -  & - & - & - \\
0 & 0 & | & C  & G_2 \\
B^* & G_2^* & | & 0 & 0 
\end{bNiceMatrix} \stackrel{P_{\text{blue}}}{\looparrowright} 
\begin{bNiceMatrix}
0 & G_1 & | & 0 & B \\
G_1^* & 0 & | & C^* & 0 \\ 
- & -  & - & - & - \\
0 &  C & |& 0  & G_2 \\
B^* & 0 & | & G_2^* & 0 
\end{bNiceMatrix} =: \Hp}
\end{equation}
In this way, the diagonal blocks are written in terms of only the corresponding diagonal blocks of $H$ and, at the same time, the off-diagonal blocks do not depend on them. Moreover, since the permutations $P_{\text{red}}$ and $P_{\text{blue}}$ are the inverse of each other ($P_{\text{red}} = P_{\text{blue}}^{-1}$), the resulting matrix  $\Hp$ is 
similar (i.e., it has the same eigenvalues) to $H_{JW}$ \cite[Cor. 1.3.4]{hornMatrixAnalysis2018}:
\begin{displaymath}
    \Hp = P_{\text{red}} H_{JW} P_{\text{red}}^{-1}.
\end{displaymath}
By the same process that led us from $H$ to $\Hp$, we can go from a perturbed version
\begin{displaymath}
    \hat{H}:= H + \begin{bmatrix}
        F_{11} & F_{12} \\ F_{21} & F_{22}
    \end{bmatrix}=: H + F,
\end{displaymath}
to the symmetric perturbed matrix
\begin{equation}
\label{eq:permutedPerturbation}
    \Hphat = \Hp + \begin{bmatrix}
0 & F_{11} & | & 0 & F_{12} \\
F_{11}^* & 0 & | & F_{21}^* & 0 \\ 
- & -  & - & - & - \\
0 &  F_{21} & | & 0 & F_{22} \\
F_{12}^* & 0 & | & F_{22}^* & 0 
\end{bmatrix}=: \Hp + \Fp .
\end{equation}
We can now apply the result \cref{eq:Yujion2x2} on the matrices $\Hp$ with eigenvalues $\{\lambda_j\}_{j=1}^{m+n}$ and $\Hphat$ with eigenvalues $\{\hat{\lambda}_j\}_{j=1}^{m+n}$. Thus, for $j=1, \dots, m+n$, we can define 
\begin{displaymath}
    \tau_j = \left( \frac{\left\Vert \begin{bmatrix}
        0 & C \\
        B^* & 0
    \end{bmatrix}\right\Vert_2 +
    \left\Vert \begin{bmatrix}
        0 & F_{21} \\
        F_{12}^* & 0
    \end{bmatrix} \right\Vert_2}
    {\displaystyle \min_\minindex|\lambda_j - \lambda_\minindex\left(\begin{bmatrix}
        0 & G_2 \\
        G_2^* & 0
    \end{bmatrix} \right)|
    - 2\left\Vert \Fp \right\Vert_2}\right),
\end{displaymath}
and conclude that, for all $j$ for which $\tau_j > 0$, we have
\begin{displaymath}
\label{eq:AdaptedYujion2x2}
    | \lambda_j - \hat{\lambda}_j | \leq \left\Vert \begin{bmatrix}
        0 & F_{11} \\
        F_{11}^* & 0
    \end{bmatrix}\right\Vert_2 + 2\left\Vert \begin{bmatrix}
        0 & F_{21} \\
        F_{12}^* & 0
    \end{bmatrix} \right\Vert_2 \tau_j +
      \left\Vert \begin{bmatrix}
        0 & F_{22} \\
        F_{22}^* & 0
    \end{bmatrix} \right\Vert_2 
    \tau_j^2.\\
\end{displaymath}
The next step is to write this bound in terms of the original blocks to find a bound for the original problem. By the Jordan-Wielandt theorem and, in particular, by \cref{eq:eigenJW}, the first $n$ eigenvalues of  $\Hp$ and $\Hphat$ are exactly the singular values $\{\sigma_i\}_{i=1}^n$ and $\{\hat{\sigma}_i\}_{i=1}^n$ of $H$ and $\hat{H}$, respectively, and, analogously, $\|\Fp\|_2 = \|F\|_2$. Then, since 
\begin{displaymath}
\left\Vert \begin{bmatrix}
        0 & M_1 \\
        M_2 & 0
    \end{bmatrix}\right\Vert_2 = \max \{\|M_1\|_2,\|M_2\|_2\},
\end{displaymath}
we obtain the desired generalization of \cite[Thm. 3.2]{nakatsukasaEigenvaluePerturbationBounds2012}:
\begin{theorem} 
\label{thm:2by2}
Consider a $2\times 2$ block matrix 
\begin{equation}
\label{def:H2by2}
    H := \begin{bmatrix}
        G_1 & B \\ C & G_2
    \end{bmatrix}
\end{equation}
with singular values $\{\sigma_i\}_{i=1}^{n}$ and the perturbed matrix
\begin{displaymath}
    \hat{H}:= H + \begin{bmatrix}
        F_{11} & F_{12} \\ F_{21} & F_{22}
    \end{bmatrix}=: H + F,
\end{displaymath}
with singular values $\{\hat{\sigma}_i\}_{i=1}^{n}$. Define, for $i=1, \dots, n$,
\begin{displaymath}
    \tau_i = \left( \frac{\max \{ \|B\|_2, \|C\|_2 \} + \max \{\left\Vert  F_{12}  \right\Vert_2, \left\Vert  F_{21}  \right\Vert_2\}}
    {\displaystyle \min_\minindex|\sigma_i - \sigma_\minindex\left(G_2 \right)|
    - 2\left\Vert  F \right\Vert_2}\right).
\end{displaymath}
Then, for all $i$ for which $\tau_i >0$, it holds
\begin{equation}
\label{eq:GeneralisedYujion2x2}
    | \sigma_i - \hat{\sigma}_i | \leq \left\Vert  F_{11}  \right\Vert_2 + 2 \max \{\left\Vert  F_{12}  \right\Vert_2, \left\Vert  F_{21}  \right\Vert_2\} \tau_i +  \left\Vert  F_{22}  \right\Vert_2 \tau_i^2.
\end{equation}
\end{theorem}
By construction, this result inherits all the properties of the symmetric case studied in \cite{nakatsukasaEigenvaluePerturbationBounds2012}. 
In particular, we need  $\tau_i < 1$ for the bound to be better than Weyl's inequality \cref{eq:Weyl}, and if $\|F_{11}\| \ll \|F\|$ and $\sigma_i$ is far from the spectrum of $G_2$ then $\tau_i \ll 1$. 
\begin{corollary}[Perturbation ``on the right''] 
\label{cor:pertOnRight}
Consider a matrix as in \cref{def:H2by2} and assume that just the two right blocks are perturbed, i.e.
\begin{displaymath}
    F = \begin{bmatrix}
        0 & F_1 \\ 0 & F_2
    \end{bmatrix}.
    \end{displaymath}
Define, for $i=1,\dots,n$,
    \begin{displaymath}
    \tau_i = \left( \frac{\max \{ \|B\|_2, \|C\|_2 \} + \left\Vert  F_1  \right\Vert_2}
    {\displaystyle \min_\minindex|\sigma_i - \sigma_\minindex\left(G_2 \right)|
    - 2\left\Vert  F \right\Vert_2}\right).
\end{displaymath}
    Then, for all $i$ for which $\tau_i>0$, \cref{eq:GeneralisedYujion2x2} trivially becomes
    \begin{equation}
    \label{eq:pertOnRight}
    | \sigma_i - \hat{\sigma}_i | \leq  2 \left\Vert  F_1  \right\Vert_2 \tau_i +  \left\Vert  F_2  \right\Vert_2 \tau_i^2.
\end{equation}
\end{corollary}
As shown in \Cref{sec:GNPerturbation}, the \GNa{} can be seen, up to orthogonal transformations, as a perturbation ``on the right'' of the correspondent original matrix. While in the general matrix perturbation theory framework, $\|F\|_2$ (i.e., Weyl's theorem) is not guaranteed to be small, in the extraction of singular values framework, if good approximate subspaces are available, $F$ will be in general quite small. Indeed, with $A$ low rank, this will follow from the fact that \GNa{} performs a low-rank approximation. However, it is the (non-zero) structure of $F$ and the spectral distribution of the (transformed) matrix $\bA$ that will determine the quality of the approximation. The denominator of $\tau_i$ will be big since the structure of the transformed matrix $\bA$ typically implies that the leading singular values are well separated from the spectrum of its $(2,2)$ block $\bA_{22}$. Thus, in \Cref{subsec:AppltoGN}, we use \Cref{cor:pertOnRight} to obtain a bound on the \GNa, and, in \Cref{sec:NumIll}, we numerically illustrate its behavior.

\section{Application to Methods for Extracting Singular Values}
\label{sec:ApplToMethods}
The discussion of \Cref{sec:MatPert} can be used to analyze methods for extracting singular values through a matrix perturbation approach. In the following, we show how this can be done for the \GNa{} and then generalize the derivation to other extracting methods.
%%%%%%%%%%%%%%%%%%%%%%%%%%%%%%%%%%%%%%%%%%%%%%%%%%%%%%%%%%%%%%%%%%%%%%%% APPLICATION TO GN %%%%%%%%%%%%%%%%%%%%%%%%%%%%%%%%%%%%%%%%%%%%%%%%%%%%%%%%%%%%%%%%%%%%%%%%
\subsection{Application to Generalized Nyström} \label{subsec:AppltoGN}
In \Cref{sec:GNPerturbation} we saw that, with the appropriate orthogonal transformations, the \GNa{} can be interpreted as a perturbation of the original matrix in the top and bottom right blocks, i.e., as a perturbation ``on the right''. Thus, we can apply \Cref{cor:pertOnRight} to the transformed matrix $\bA$ and its \GNa{} with sampling matrices \cref{eq:GNIdentities-final}, obtaining, by \cref{eq:transfSingVal_eq_SingVal}, a bound on the difference between the extracted singular values $\sigma^{\mtiny{GN}}_i:=\sigma_i(A_{\mtiny{GN,\tV,\tU}})$ and the exact singular values $\sigma_i:=  \sigma_i(A)$, for $i=1, \dots, n$. 

\begin{corollary}
\label{thm:GNbound}
    Consider an $m\times n$ matrix $A$ with singular values $\{\sigma_i\}_{i=1}^n$ and approximations $\tU$ and $\tV$ to its leading singular subspaces. Consider the \GNa\ 
    $A_{\mtiny{GN,\tV,\tU}} = A\tV (\tU^*A\tV)^\dagger \tU^* A$ and its singular values $\{\sigma^{\mtiny{GN}}_i\}_{i=1}^n$. Define $\bA = Q_1^*AQ_2$, with $Q_{1}= \begin{bmatrix} \tU & \tU_{\perp}\end{bmatrix}$ and $Q_{2}= \begin{bmatrix} \tV & \tV_{\perp}\end{bmatrix}$, and its blocks as in \eqref{def:Abar}. Let, for $i=1, \dots, n$,
    \begin{displaymath}
    \tau_i = \frac{\max \{ \|\bA_{12}\|_2, \|\bA_{21}\|_2 \} + \left\Vert  \bA_{12} - \bA_{11}\bA_{11}^\dag \bA_{12} \right\Vert_2} {\displaystyle \min_\minindex|\sigma_i - \sigma_\minindex\left(\bA_{22} \right)|- 2\left\Vert  E_{\mtiny{GN}} \right\Vert_2}.
\end{displaymath}
Then, for all $i$ for which $\tau_i>0$, it holds
\begin{equation}
\label{eq:GNbound}
    | \sigma_i - \sigma^{\mtiny{GN}}_i | \leq  2\left\Vert  \bA_{12} - \bA_{11}\bA_{11}^\dag \bA_{12}  \right\Vert_2 \tau_i + \left\Vert  \bA_{22} - \bA_{21}\bA_{11}^\dag \bA_{12} \right\Vert_2 \tau_i^2.
\end{equation}
\end{corollary}
As already mentioned, in the case where $\tU$ and $\tV$ are given with the same column size (i.e., without oversampling, $\ell =0$), we have $\| \bA_{12} - \bA_{11}\bA_{11}^\dag \bA_{12}  \|_2 = 0 $. Thus, the bound can be slightly simplified. However, we note once again that, since the block sizes change compared to the oversampled case, this does not mean that $\ell =0$ implies better approximations.

\subsection{Application to Other Extraction Methods}
\label{subsec:Appltoother}
The \GNa{} is not the only method that can be analyzed through a matrix perturbation approach. Indeed, all methods mentioned in the introduction can be interpreted as a perturbation of the original matrix. Moreover, also in these other cases, we can achieve, by orthogonal transformations, a special structure of the perturbation that can be exploited, using \Cref{thm:2by2}, to derive a bound on the singular values extraction error. Therefore, as explored in \Cref{subsec:Comparison}, our analysis can be seen as a tool for potential comparisons between the different strategies. For example, let us consider the \RR{} method  
    \begin{displaymath}
      \sigma_i := \sigma_i(A) \approx \sigma_i(\tU^*A\tV) =:\sigma_i^{\mtiny{RR}}, \quad i=1, \dots ,r.
    \end{displaymath}
    Then, considering the orthogonal transformations
       \begin{displaymath} 
       Q_{1}= \begin{bmatrix} \tU & \tU_{\perp}\end{bmatrix}, \quad Q_{2}= \begin{bmatrix} \tV & \tV_{\perp}\end{bmatrix},
 \end{displaymath}
    we have
    \begin{displaymath}
    (Q_1^*AQ_2)_{\mtiny{ RR,\begin{bmatrix}I_r \\0\end{bmatrix},\begin{bmatrix}I_{r + \ell} \\ 0\end{bmatrix}}} = \begin{bmatrix}
        I_{r + \ell} & 0
    \end{bmatrix}Q_1^*AQ_2\begin{bmatrix}
        I_{r} \\ 0
    \end{bmatrix} = \tU^*A\tV = A_{\mtiny{ RR,\tV,\tU}}.
    \end{displaymath}
This means that instead of studying the original extraction error, we can study the error of the \RR{} method with sampling matrices applied to $\bA = Q_1^*AQ_2$, i.e., for $i=1,\dots, r$,
\begin{displaymath}
\begin{aligned}
            |\sigma_i - \sigma_i^{\mtiny{RR}}| &= |\sigma_i(\bA) - \sigma_i(\bA_{\mtiny{ RR,\begin{bmatrix}I_r \\0\end{bmatrix},\begin{bmatrix}I_{r + \ell} \\ 0\end{bmatrix}}})|  = |\sigma_i(\bA) - \sigma_i(\bA_{11})| \\
            &=  |\sigma_i(\bA) - \sigma_i\left(\begin{bmatrix}
                \bA_{11} & 0 \\ 0 & 0
            \end{bmatrix}\right)| = |\sigma_i(\bA) - \sigma_i\biggl(\bA - \underbrace{\begin{bmatrix}
                0 & \bA_{12} \\ \bA_{21} & \bA_{22}
            \end{bmatrix}}_{=: E_{\mtiny{RR}}} \biggr)|.
\end{aligned}
    \end{displaymath}
Therefore, it is possible to apply \Cref{thm:2by2} and, as a consequence, obtain that, for all $i$ for which $\tau_i^{\mtiny{ RR}}:= \frac{2\max \{ \|\bA_{12}\|_2, \|\bA_{21}\|_2 \}}{(\min_\minindex|\sigma_i - \sigma_\minindex\left(\bA_{22} \right)|- 2\left\Vert  E_{\mtiny{RR}} \right\Vert_2)} > 0$,
\begin{equation}
\label{eq:RRbound}
\begin{aligned}
     |\sigma_i - \sigma_i^{\mtiny{RR}}| \leq \, &4 \frac{\max \{ \|\bA_{12}\|_2, \|\bA_{21}\|_2 \}^2}{\displaystyle \min_\minindex|\sigma_i - \sigma_\minindex\left(\bA_{22} \right)|- 2\left\Vert  E_{\mtiny{RR}} \right\Vert_2}\\
     &+ \left\Vert  \bA_{22}  \right\Vert_2\frac{4\max \{ \|\bA_{12}\|_2, \|\bA_{21}\|_2 \}^2}{(\displaystyle \min_\minindex|\sigma_i - \sigma_\minindex\left(\bA_{22} \right)|- 2\left\Vert  E_{\mtiny{RR}} \right\Vert_2)^2}.
\end{aligned}
\end{equation}
Note that the \RR{} bound is largely influenced by the non-zero off-diagonal perturbations and that the significant difference in accuracy between \GN{} and \RR{} may arise from the fact that, unlike \RR, off-diagonal perturbations in \GN{} can be much smaller (or zero if oversampling is not used).

To obtain a bound for the SVD approximation we note that
\begin{displaymath}
\begin{aligned}
|\sigma_i - \sigma_i^{\mtiny{SVD}} | &= |\sigma_i - \sigma_i(A\tV ) | = |\sigma_i(AQ_2) - \sigma_i((AQ_2)_{\mtiny{ SVD,\begin{bmatrix}I_r \\0\end{bmatrix}}}) |\\ 
&\stackrel{\mathclap{\tA := AQ_2}}{=} \quad|\sigma_i(\tA) - \sigma_i(\tA - \underbrace{\begin{bmatrix}
    0 & \tA_2
\end{bmatrix}}_{=:E_{\mtiny{SVD}}})|,
\end{aligned}
\end{displaymath}
where $\tA_2$ is the $m \times (n-r)$ right block of $\tA$. Then, using \Cref{thm:2by2}, we obtain, for all $i$ for which $\tau_i^{\mtiny{ SVD}}:=\frac{2\|\tA_2\|_2}{\sigma_i - 2\|E_{\mtiny{SVD}}\|_2}> 0$,
\begin{displaymath}
    |\sigma_i - \sigma_i^{\mtiny{SVD}} | \leq 4 \frac{\|\tA_2\|^2_2}{\sigma_i - 2\|E_{\mtiny{SVD}}\|_2}.
\end{displaymath}
Finally, to obtain a bound for the \HMTa, we notice that, as shown in \cite[Table 1]{nakatsukasaFastStableRandomized2020}, if $Q$ is the orthogonal matrix from the QR factorization of $A\tV$, then:
\begin{equation}
\label{eq:HMTasGN}
    A_{\mtiny{HMT}} := QQ^*A = A_{\mtiny{ GN, \tV, A\tV}}.
\end{equation}
Thus, the \HMTa{} can be seen as a specific \GNa{} enabling us to use \Cref{thm:GNbound}. Gu~\cite{gu2015subspace} and Saibaba~\cite{saibaba2019randomized} provide bounds for the extraction of singular values by the HMT method. In \cref{subsec:Saibaba}, we show that, in contrast with bound \eqref{eq:GNbound}, while these bounds are particularly sharp when $\tV$ is a random Gaussian matrix, they are less precise when $\tV$ has information of $A$.
%%%%%%%%%%%%%%%%%%%%%%%%%%%%%%%%%%%%%%%%%%%%%%%%%%%%%%%%%%%%%%%%%%%%%%%%%%%%%%%%%%%%%%%%%%%%%%%%%%%%%%%%%%%%%%%%%%%%%%%%%%%%%%%%%%%%%%%%%%%%%%%%
%%%%%%%%%%%%%%%%%%%%%%%%%%%%%%%%%%%%%%%%%%%%%%%%%%%%%%%%%%%%% NUMERICAL ILLUSTRATIONS %%%%%%%%%%%%%%%%%%%%%%%%%%%%%%%%%%%%%%%%%%%%%%%%%%%%%%%%%%
%%%%%%%%%%%%%%%%%%%%%%%%%%%%%%%%%%%%%%%%%%%%%%%%%%%%%%%%%%%%%%%%%%%%%%%%%%%%%%%%%%%%%%%%%%%%%%%%%%%%%%%%%%%%%%%%%%%%%%%%%%%%%%%%%%%%%%%%%%%%%%%% 
\section{Numerical Illustrations and Comparisons} \label{sec:NumIll}
\begin{algorithm}[t]
    \tt 
    \small
    \caption{The \GNa{} (\textsc{Matlab} notation): given $A\in \mathbb{R}^{m \times n}, \tV \in \mathbb{R}^{n\times r}$ and $\tU \in \mathbb{R}^{m\times (r + \ell)}$, approximate the first $r$ singular values of $A$.}
\begin{algorithmic}[1]
\State AtV = A * tV;
\State tUA = tU' * A;
\State [$\sim$, R1] = qr(AtV, 0); \Comment{QR factorization of external factor}
\State [$\sim$, R2] = qr(tUA', 0); \Comment{QR factorization of external factor}
\State [Q3, R3] = qr(tU' * AtV, 0);                 \Comment{QR factorization to compute pseudoinverse}
\State P = Q3' * R2';        \Comment{Note: The order of operations matters!}
\State AGN = (R1/R3) * P;            \Comment{Compute \GNa}
\State svGN = svd(AGN);                                      \Comment{Use build-in function to compute singular values}
\end{algorithmic}
\label{algo:GN}
\end{algorithm}
Having derived results that respect the original goal of adaptability depending on each singular value and the use of the perturbation structure, we now numerically illustrate the sharpness of the derived bounds, compare different methods, and explore the possibility of making the derived bounds computable in practice.

We construct a $1000\times 1000$ matrix $A$ starting from its SVD $A = U \Sigma V^*$, with $U$ and $V$ Haar matrices~\cite{meckes2019random}. We construct the singular subspace approximations from Gaussian matrices subjected to one power iteration. That is, taking Gaussian matrices $\Omega_1$ and $\Omega_2$, we define $\tV \in \mathbb{R}^{1000\times 200}$ and $\tU \in \mathbb{R}^{1000\times (200+\ell)}$ as the orthonormal matrices of the QR factorization of $A^*\Omega_1$ and $A\Omega_2$, respectively. There are various alternative strategies to generate approximate singular subspaces, including subspace iteration and Krylov methods. This means, e.g., generating $\tilde{V}$ and $\tilde{U}$ as the orthogonal basis computed by the Lanczos algorithm on $A^*A$ and $AA^*$ after $r$ and $r+\ell$ iterations, respectivetly. We have also experimented with the latter, but ultimately considered the case described above which is of more interest.
Indeed, when generating subspaces by Krylov methods, $A$ is applied repetitively, thus $AA^*v$ is contained in the span of $V$ for all but the last vector in $[v_0,AA^*v_0,\ldots, (AA^*)^rv]$. Thus, the second multiplication by $A$, typical of \GN{}, has less impact than in other situations, leading to rather small improvements on other single pass methods.

We consider two singular value distributions. Firstly, the singular values will follow an exponential decay:
    \begin{equation}
    \label{eq:expDecay}
        \sigma_i = e^{-(i-1)\frac{30}{999}\ln{10}},
    \end{equation}
    where the largest singular value is $1$, the smallest $10^{-30}$, and $\sigma_r = \sigma_{200}$ is of the order of $10^{-6}$. This can be implemented by the \textsc{Matlab} function \verb|logspace(-30,0,n)'|.\footnote{All the implementations are available in the Supplementary Material.} Secondly, the singular values will follow an algebraic decay:
    \begin{displaymath}
        \sigma_i = \left( \frac{1}{i} \right)^4,
    \end{displaymath}
    where the largest singular value is $1$, the smallest $10^{-12}$, and $\sigma_r = \sigma_{200}= 6.25 \cdot 10^{-10}$. 
    
    We compute the \GNa{} and the extracted singular values by \Cref{algo:GN}, in which we exploit the structure of the approximation instead of forming the $m\times n$ approximate matrix. Namely, first, we compute the QR (LQ) factorization of the external factors (lines 3 and 4)
    \begin{displaymath}
    \begin{aligned}
        & A\tV = Q_1 R_1 \\ 
        & \tU A^{*} = Q_2R_2 \implies \tU^{*}A = R_2^{*}Q_2^{*}.
    \end{aligned}
    \end{displaymath}
    To compute the pseudoinverse in \cref{def:GN}, we use the QR factorization of the matrix to invert $\tU^{*}A\tV = Q_3R_3$ (line 5).
    Then, since singular values do not change under orthogonal transformations, we extract the positive ones by looking only at the small matrix formed by the R-factors and the pseudoinverse of the core matrix (lines 6 and 7)
    \begin{displaymath}
    \begin{aligned}
        \sigma^{\mtiny{GN}} &= \sigma(Q_1R_1 (Q_3R_3)^\dagger R_2^{*} Q_2^{*}) \\
        & = \sigma(Q_1R_1 R_3^\dagger Q_3^{*} R_2^{*} Q_2^{*}) \\
        & = \sigma(R_1 R_3^\dagger Q_3^{*} R_2^{*}).
    \end{aligned}
    \end{displaymath}
    The singular values of the small matrix are then computed by the  \textsc{Matlab} function \verb|svd|. Alternatively, it is possible to use the Jacobi SVD from \cite{drmavc2008newI,drmavc2008newII}, which can compute singular values with high relative accuracy for certain class of matrices.
    In the whole procedure, we need to perform two matrix-matrix multiplications by $A$ (lines 1 and 2) but, since they are not sequential (and hence parallelizable), we need to access $A$ only once. 
    
    To understand the complexity of this procedure, we start by considering the cost $\Nr$ and $N_{r+\ell}$  of computing $A\tV$ or $\tU^*A$, respectively. This varies depending on the properties of $A$, e.g. sparsity, and those of the approximate singular subspaces. The QR factorization (and consequently computing the SVD) of an $m\times n$ ($m>n$) matrix costs $\calO(mn^{2})$. Therefore, from \Cref{algo:GN}, we see that the cost of the \GN{} procedure is dominated by the matrix multiplications with $A$ of cost $N_{2r+\ell}:=\Nr + N_{r+\ell}$, and by the QR factorizations needed to extract singular values, particularly the one in line 3. Hence, we obtain an overall cost of $N_{2r+\ell} +\calO((m+n)r^2)$. 
    \subsection{Equal Column Size Singular Approximate Subspace \texorpdfstring{($\boldsymbol{\ell = 0}$)}{}} 
    \label{subsec:No-Oversample} 
            \begin{figure}[t]
    \begin{subfigure}{.5\textwidth}
          \centering
          \includegraphics[width=\linewidth]{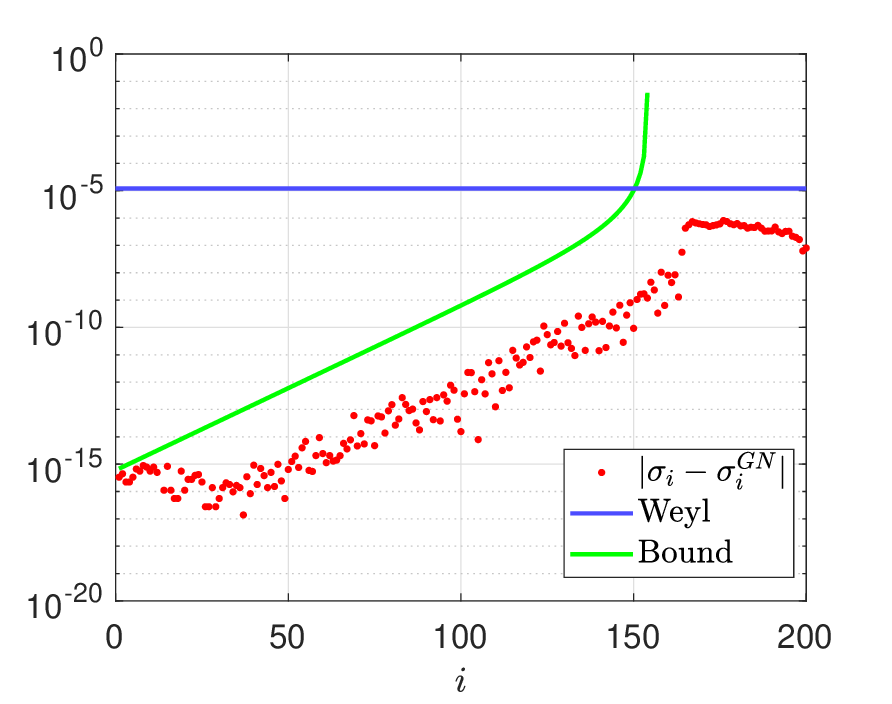}
          \caption{Exponential Decay}
          \label{fig:GN_exp}
    \end{subfigure}%
    \begin{subfigure}{.5\textwidth}
          \centering
          \includegraphics[width=\linewidth]{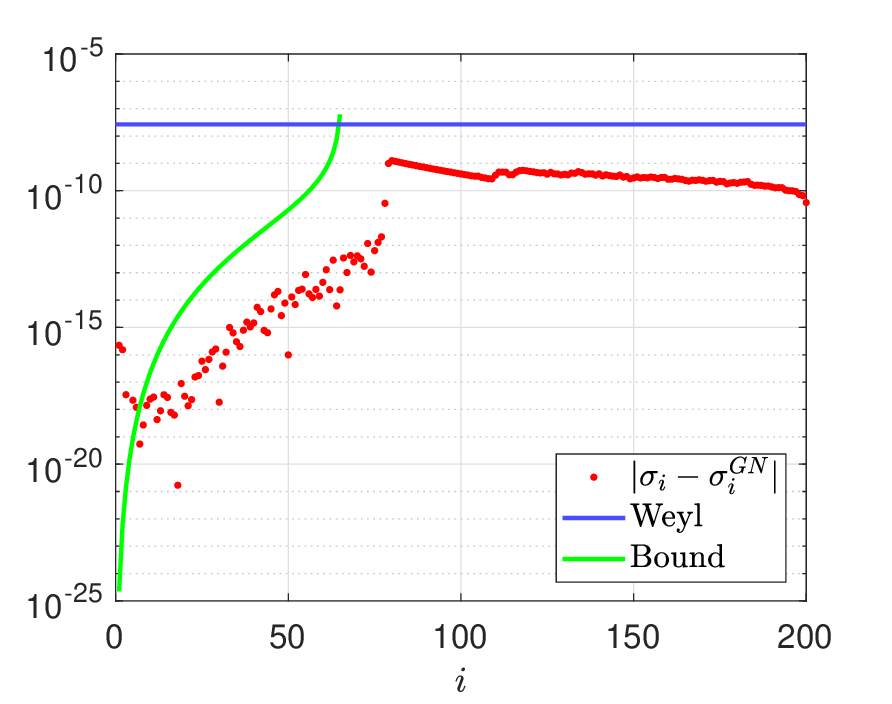}
          \caption{Algebraic Decay}
          \label{fig:GN_alg}
    \end{subfigure}
    \caption{\textbf{Singular value approximation errors of \GN{} without oversampling.} We show the error in extracting the first $r = 200$ singular values by \GNa{} with $\ell=0$ (red dots), i.e., $|\sigma_i - \sigma_i^{\mtiny{GN}}|$ for $i=1,\dots, 200$, Weyl's inequality \cref{eq:Weyl} (blue), and bound \cref{eq:GNbound} (green). Both exponential (a) and algebraic (b) singular values decays are considered.}
    \label{fig:GN}
    \end{figure}
    Let us consider the case without oversampling, i.e., $\ell = 0$, and, where $\tU$ and $\tV$ belong to $\mathbb{R}^{1000\times 200}$. In \Cref{fig:GN}, we observe that, for both exponential and algebraic decay, the bound in \Cref{thm:GNbound} is significantly better than Weyl's inequality and that it can better predict the slope of the error. It is important to note that the apparent misbehavior of the bound for the very first singular value errors is simply due to the bounds falling below machine precision; as our analysis does not account for roundoff errors.
    \subsection{Different Column Size Approximate Singular Subspace \texorpdfstring{($\boldsymbol{\ell = 0.5r}$)}{}} \label{subsec:Oversample}
     \begin{figure}
    \begin{subfigure}{.5\textwidth}
          \centering
          \includegraphics[width=\linewidth]{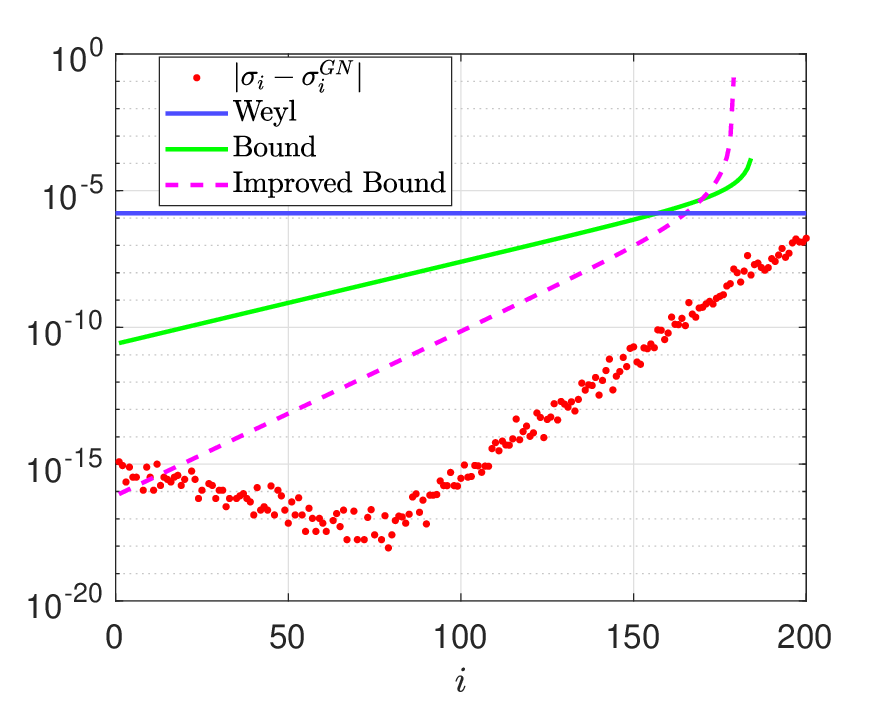}
          \caption{Exponential Decay}
          \label{fig:GN_exp_ov}
    \end{subfigure}%
    \begin{subfigure}{.5\textwidth}
          \centering
          \includegraphics[width=\linewidth]{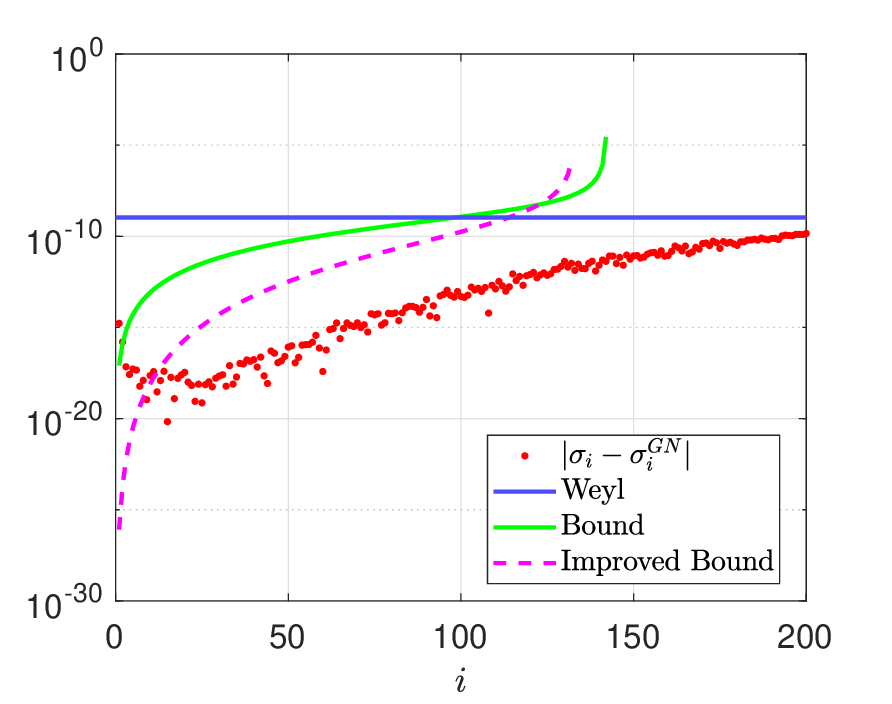}
          \caption{Algebraic Decay}
          \label{fig:GN_alg_ov}
    \end{subfigure}
    \caption{\textbf{Singular value approximation errors of \GN{} with oversampling.} We show the error in extracting the first $r = 200$ singular values by \GNa{} with $\ell=1.5r$ (red dots), i.e., $|\sigma_i - \sigma_i^{\mtiny{GN}}|$ for $i=1,\dots, 200$, Weyl's inequality \cref{eq:Weyl} (blue), bound \cref{eq:GNbound} (green), and the improved bound (magenta). Both exponential (a) and algebraic (b) singular values decays are considered.}
    \label{fig:GN_ov}
    \end{figure}
    Let us consider the case in which the approximations are given with different column sizes. In particular, we set $r+\ell = 1.5r$ following the recommendation in \cite{nakatsukasaFastStableRandomized2020}, so that, for $r=200$, we take $\tV \in \mathbb{R}^{1000\times 200}$ and $\tU \in \mathbb{R}^{1000\times 300}$. \Cref{fig:GN_ov} shows that for both algebraic and exponential singular value decays, bound \cref{eq:GNbound} is still significantly better than Weyl's. However, the sharpness observed in the previous experiment seems partially lost. Qualitatively speaking, this could be caused by the need to compute the pseudoinverse of a rectangular matrix ($\bA_{11}^\dagger$) and by not achieving the typical structure discussed in \Cref{sec:GNPerturbation} when constructing $\bA$. In view of this, we propose a heuristics-based technique that, after extensive numerical tests, seems to improve this bound and usually estimate the error reasonably accurately.
    
   \paragraph{Heuristic-based improvement}      
   The goal is to reduce ourselves to a situation similar to  the non-oversampling framework. This means that we aim at modifying the sizes of the blocks used to compute bound \eqref{eq:GNbound}, so as to have a square $(1,1)$ block. Moreover, we rely on the fact that oversampling generally improves \GN{} approximations. Naively modifying the sizes of $\bA$ blocks (as in \cref{def:Abar}) would probably ruin even more the desired structure described in \Cref{sec:GNPerturbation}. Thus, we first apply orthogonal transformations to the (already transformed) matrix $\bA$ to diagonalize its $(1,1)$ block, $\bA_{11} = X\Sigma Y^{*}$, and obtain the matrix \begin{displaymath}
        \bbA := \begin{bmatrix}
            X^{*} & 0 \\ 0 & I
        \end{bmatrix} \bA \begin{bmatrix}
            Y & 0 \\ 0 & I
        \end{bmatrix}.
    \end{displaymath}
    Since $\bbA$ is an orthogonal transformation of $\bA$, for the same reasons illustrated at the start of \Cref{sec:GNPerturbation}, we have
    \begin{displaymath}
        | \sigma_i(\bA) - \sigma_i(\bA_{\mtiny{ GN,\begin{bmatrix}I_r \\0\end{bmatrix},\begin{bmatrix}I_{r + \ell} \\ 0\end{bmatrix}}})| = | \sigma_i(\bbA) - \sigma_i(\bbA_{\mtiny{ GN,\begin{bmatrix}Y^{*} \\0\end{bmatrix},\begin{bmatrix}X^{*} \\ 0\end{bmatrix}}}) |.
    \end{displaymath}
    Moreover, due to the orthogonal invariance of the $2$-norm, bound \cref{eq:GNbound} reduces to the same expression when computed for $| \sigma_i(\bbA) - \sigma_i(\bbA_{\mtiny{ GN,\begin{bmatrix}Y^{*} \\0\end{bmatrix},\begin{bmatrix}X^{*} \\ 0\end{bmatrix}}}) |$ and for $| \sigma_i(\bbA) - \sigma_i(\bbA_{\mtiny{ GN,\begin{bmatrix}I_r \\0\end{bmatrix},\begin{bmatrix}I_{r + \ell} \\ 0\end{bmatrix}}})|$. Finally, note that the latter is, in general, smaller than the error obtained if no oversampling is employed, i.e., $| \sigma_i(\bbA) - \sigma_i(\bbA_{\mtiny{ GN,\begin{bmatrix}I_r \\0\end{bmatrix},\begin{bmatrix}I_{r} \\ 0\end{bmatrix}}})|$. Thus, we can expect that bound \cref{eq:GNbound} computed with blocks of $\bbA$ whose sizes are chosen so that $\bbA_{11}$ is $r\times r$ (in the specific example $200\times 200$) is going to bound the original accuracy error.

\begin{figure}
       \centering
       \includegraphics[width=1\linewidth]{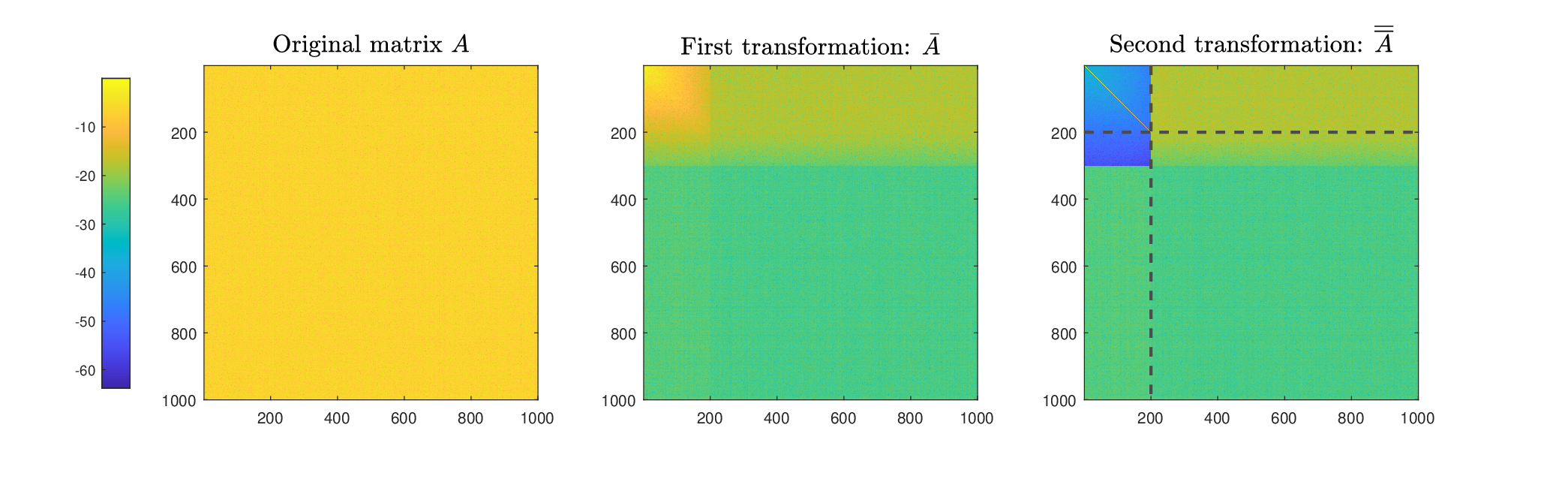}
       \caption{\textbf{Graphic representation of the process to obtain an heuristic-based bound improvement.} The figures show the magnitude of the entries in logaritmic scale of the matrices $A,\bA$, and $\bbA$, with $A$ generated with exponentially decaying singular values. The dotted lines in the figure corresponding to $\bbA$ indicate the choice of the blocks sizes used in computing the heuristic-based bound.}
       \label{fig:heuristic}
   \end{figure}

    In \Cref{fig:heuristic}, we illustrate the effects that this process has on the matrix with exponentially decaying singular values taken into account above. \Cref{fig:GN_ov} shows the behavior of bound \cref{eq:GNbound} and the heuristically improved version for both exponential and algebraic singular value decays. Note that, as shown by both \Cref{fig:GN} and \cref{fig:GN_ov}, in general, there exists an $i < r$ for which the presented bounds became worse than Weyl's one. This is because, when approaching the $r$-th singular value, the gap becomes smaller, making $\tau_i$ greater than $1$. Thus, we suggest considering the bound given by the minimum, for each $i$, between \cref{eq:GNbound} and Weyl's inequality.
    
%%%%%%%%%%%%%%%%%%%%%%%%%%%%%%%%%%%%%%%%%%%%%%%%%%%%%%%%%%%%%%%%%%%%%%%% METHODS COMPARISON %%%%%%%%%%%%%%%%%%%%%%%%%%%%%%%%%%%%%%%%%%%%%%%%%%%%%%%%%%%%%%%%%%%%%%%%
    \subsection{Methods Comparison} \label{subsec:Comparison}
We implement and compare the various methods. We compute the \RR{} and the SVD approximations by \Cref{algo:RR} and \Cref{algo:SVD}, respectively. In both cases, we first construct the approximate matrix $A_{\mtiny{RR}}$ or $A_{\mtiny{SVD}}$ and then use the \textsc{Matlab} build-in function \verb|svd| to compute the singular values. We stress that, in both cases, we need to access the matrix $A$ only once, i.e., these are single-pass algorithms, and we need to perform only one matrix-matrix multiplication by $A$. The cost of the \RR{} procedure is given by $N_{r}+O(mr(r+\ell))$ operations plus the cost of computing the singular values of an $(r+\ell)\times r$ matrix that, assuming $\ell = \calO (r)$, is $\calO(r^3)$. While, to obtain the SVD approximation we need to compute $A\tV$, of cost $\Nr$, and compute singular values of a $m\times r$ matrix, of cost $\calO(mr^2)$. 

To implement the HMT method, \Cref{algo:HMT}, we need to access and multiply by $A$ twice, see lines 1 and 3. Moreover, the multiplication of computing $Q^*A$ almost invariably involves multiplying an unstructured, dense matrix $Q^*$ to $A$, in contrast to computing $A\tV$ or $\tU^*A$, which can be faster for structured (e.g. sparse) $\tV,\tU$. We thus distinguish the cost of computing $\tU^*A$ (which is $N_r$ or $N_{r+\ell}$) from that of $Q^*A$, which we denote by $\tNr$. 
   \begin{algorithm}[H]
    \tt 
    \small
    \caption{The \RR{} approximation (\textsc{Matlab} notation): given $A\in \mathbb{R}^{m \times n}, \tV \in \mathbb{R}^{n\times r}$ and $\tU \in \mathbb{R}^{m\times (r + \ell)}$, approximate the first $r$ singular values of $A$.}
\begin{algorithmic}[1]
\State AtV = A * tV;
\State ARR = tU' * AtV;            \Comment{Compute \RR{} approximation}
\State svRR = svd(ARR);                                      \Comment{Use build-in function to compute singular values}
\end{algorithmic}
\label{algo:RR}
\end{algorithm}
\begin{algorithm}[H]
    \tt 
    \small
    \caption{The SVD approximation (\textsc{Matlab} notation): given $A\in \mathbb{R}^{m \times n}$ and $ \tV \in \mathbb{R}^{n\times r}$, approximate the first $r$ singular values of $A$.}
\begin{algorithmic}[1]
\State ASVD = A * tV;                                          \Comment{Compute SVD approximation}
\State svSVD = svd(ASVD);                                      \Comment{Use build-in function to compute singular values}
\end{algorithmic}
\label{algo:SVD}
\end{algorithm}
\begin{algorithm}[H]
    \tt 
    \small
    \caption{The HMT method (\textsc{Matlab} notation): given $A\in \mathbb{R}^{m \times n},$ and $\tV \in \mathbb{R}^{n\times r}$, approximate the first $r$ singular values of $A$.}
\begin{algorithmic}[1]
\State AtV = A * tV; \label{algstate:HMTmultA1}
\State [Q,$\sim$] = qr(AtV, 0);           
\State AHMT = Q' * A; \label{algstate:HMTmultA2}                                      \Comment{Compute HMT approximation}
\State svHMT = svd(AHMT);                                      \Comment{Use build-in function to compute singular values}
\end{algorithmic}
\label{algo:HMT}
\end{algorithm}

\begin{figure}[t]
    \begin{subfigure}{.5\textwidth}
          \centering
          \includegraphics[width=\linewidth]{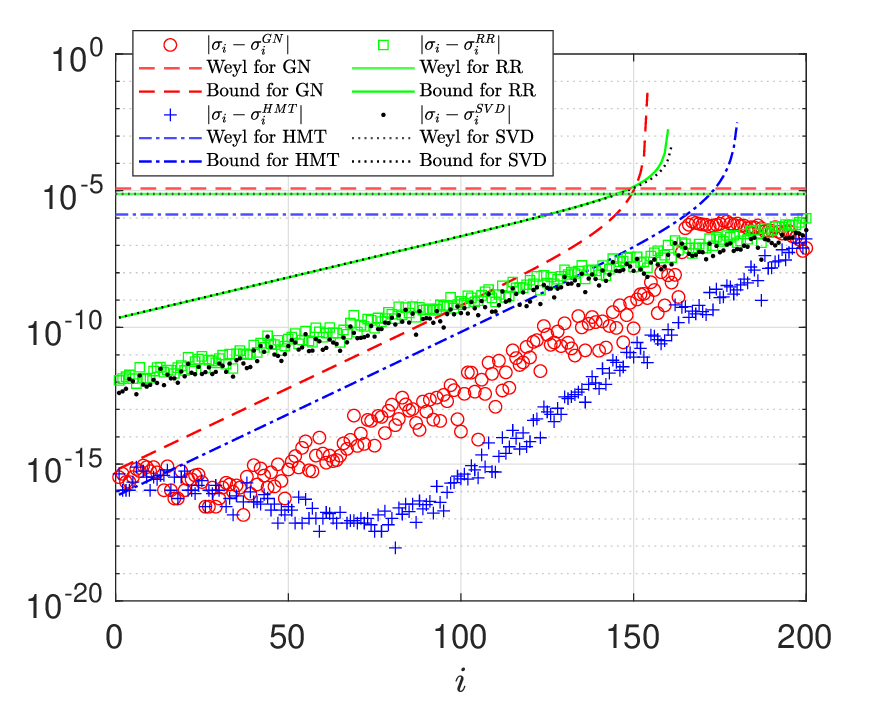}
          \caption{Without oversampling}
          \label{fig:Comparison_exp}
    \end{subfigure}%
    \begin{subfigure}{.5\textwidth}
          \centering
          \includegraphics[width=\linewidth]{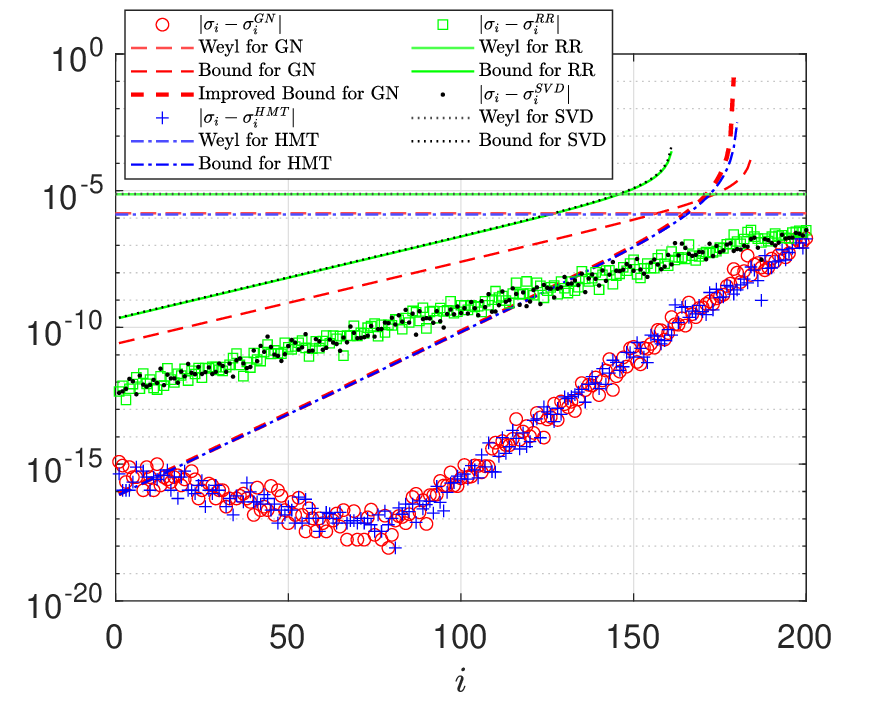}
          \caption{With oversampling}
          \label{fig:Comparison_exp_ov}
    \end{subfigure}
    \caption{\textbf{Method Comparison for Exponentially decaying exact singular values.} We show the approximation errors and relative bounds for: \GN{} (red), \RR{} (green), HMT (blue), and SVD (black). We consider both the case without (a) and with (b) oversampling.}
    \label{fig:Comparison}
    \end{figure}

In \Cref{tab:Comparison} we summarize the main features of the presented algorithms. 
\begin{table}
\centering 
\begin{tabular}{c|c|c|c}
     &  Pass & MatMul by $A$ & Complexity \\ \hline
     \RR &  Single & $1$ & $N_{r}+O(mr^2)$
%$N_{2r+\ell} +\calO( r^3)$
\\
     SVD approx. &  Single & $1$ &  $\Nr + \calO(mr^2)$\\
     Generalized Nyström &  Single & $2$ & $N_{2r+\ell} + \calO((m+n)r^2)$\\
     HMT &  Double & $2$ & $\Nr + \tNr+ \calO(mr^2)$
\end{tabular}
\caption{\textbf{Comparison of algorithms for the extraction of singular values.} $\Nr$ and $N_{r+\ell}$ are the cost of forming the product $A\tV$ or $\tU^*A$, respectively, and $N_{2r+\ell}$ is their sum. $\tNr$ is the cost of computing $Q^*A$ in the HMT method (\Cref{algo:HMT}, line 3). Here we assume $\ell = \calO(r)$. Recall that for \GN, we need not assume $\tU,\tV$ are orthonormal.
}
\label{tab:Comparison}
\end{table}
In \Cref{fig:Comparison}, we compare the methods, their extraction errors, and the corresponding bounds (derived in \Cref{subsec:Appltoother}) applied to the case where $A$ has exponentially decaying exact singular values as in \cref{eq:expDecay}.
Analogous experiments with similar results for algebraically decaying exact singular values are in Supplementary Material. Since the derived bounds seem to be good predictors of extraction errors, we can use them to explain the behavior of the methods. In particular, we consider the \GNa{} and the \RR{} method without oversampling. First of all, we note that the bounds of both methods are computed with the blocks of the same transformed matrix $\bA$. Moreover, \Cref{fig:Comparison} confirms that Weyl's inequality is not descriptive of the difference between these methods and that the perturbation norms are close to each other, i.e., $\|E_{\mtiny{GN}}\| \approx \|E_{\mtiny{RR}}\|$. Therefore, with $\ell = 0$, we have 
\begin{equation}
\label{eq:tauGNvsRR}
\begin{aligned}
 \tau_{i}^{\mtiny{GN}} &= \frac{\max \{ \|\bar{A}_{12}\|_2, \|\bar{A}_{21}\|_2 \}} {\displaystyle \min_{k}|\sigma_i - \sigma_{k}\left(\bar{A}_{22} \right)|- 2\left\Vert  E_{\mtiny{GN}} \right\Vert_{2}} \\
 &\approx \frac{1}{2}\frac{2\max \{ \|\bar{A}_{12}\|_2, \|\bar{A}_{21}\|_2 \}} {\displaystyle \min_{k}|\sigma_i - \sigma_{k}\left(\bar{A}_{22} \right)|- 2\left\Vert E_{\mtiny{RR}} \right\Vert_{2}} = \frac{1}{2}\tau_{i}^{\mtiny{RR}}.
\end{aligned}
 \end{equation}
 Now, we note that in bound \cref{eq:RRbound}, the first term is dominant, i.e.
 \begin{equation}
 \label{eq:approxRRbound}
     |\sigma_i - \sigma_i^{\mtiny{RR}}| \lesssim  2\max \{ \|\bA_{12}\|_2, \|\bA_{21}\|_2 \} \tau_i^{\mtiny{RR}},
\end{equation}
while, by using \cref{eq:tauGNvsRR} in \cref{eq:GNbound} with $\ell = 0$,
\begin{equation}
\label{eq:approxGNbound}
    |\sigma_i - \sigma_i^{\mtiny{GN}}| \lesssim \frac{1}{4}\left\Vert \bA_{22} - \bA_{21}\bA_{11}^\dag \bA_{12}\right\Vert_2 (\tau_i^{\mtiny{RR}})^2.
\end{equation}
Recall that the derived bounds intersect the corresponding Weyl bound when $\tau_i = 1$. Therefore, $\max \{ \|\bA_{12}\|_2, \|\bA_{21}\|_2 \}$ and $\| \bA_{22} - \bA_{21}\bA_{11}^\dag \bA_{12}\|_2$ are close to the respective perturbation norms and, consequently, close to each other. Thus, the presence of $(\tau_i^{\mtiny{RR}})^2$ in \cref{eq:approxGNbound} (as opposed to $\tau_i^{\mtiny{RR}}$ in \cref{eq:approxRRbound}) gives us a possible explanation of the significantly better behavior of the \GNa. 

\subsection{Comparison with bound in \cite{saibaba2019randomized}}
\label{subsec:Saibaba}
As previously mentioned, Saibaba~\cite{saibaba2019randomized} provides a bound on the accuracy of the HMT method in extracting singular values, that, given $\gamma_{j}= \frac{\sigma_{r+1}}{\sigma_{j}}$ and $$ \begin{bmatrix}\Omega_{1}\\ \Omega_2\end{bmatrix}:= V^*\tilde{V},$$ can be written as
$$\frac{|\sigma_{j}- \sigma_{j}^{HMT}|}{|\sigma_{j}|} \leq 1 - \frac{1}{\sqrt{1 + \left(\frac{\sigma_{r+1}}{\sigma_{j}}\right)^{2}\|\Omega_{2}\Omega_1^{\dagger}\|_{2}^{2}}}.$$ 
Since this is a bound on the relative error, as opposed to our bound that is on the absolute error, we divide \cref{thm:GNbound} by $|\sigma_i|$, to obtain a comparable bound.
We perform experiments with the matrix with algebraically decaying singular values as described above and no oversampling, i.e., $\ell =0$. We consider the approximate subspace $\tV$ to be, at first, as in the previous experiments, and then, randomly generated following the Haar distribution, i.e., as the orthogonal factor of a Gaussian matrix.

\begin{figure}
    \begin{subfigure}[t]{.5\textwidth}
          \centering
          \includegraphics[width=\linewidth]{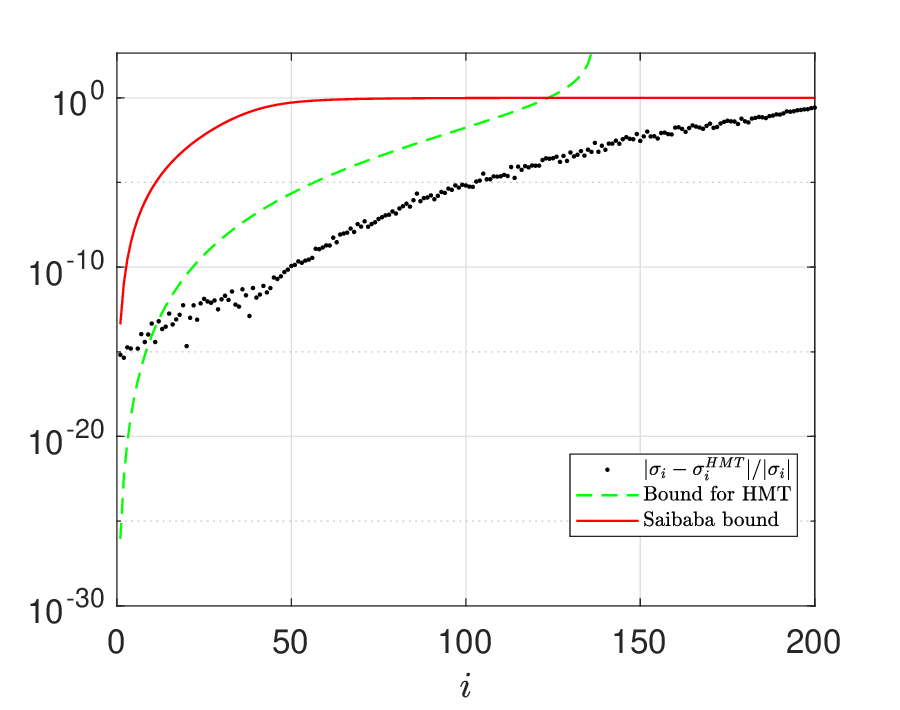}
          \caption{$\tilde{V}$ from Gaussian subject to one power iteration}
    \end{subfigure}%
    \begin{subfigure}[t]{.5\textwidth}
          \centering
          \includegraphics[width=\linewidth]{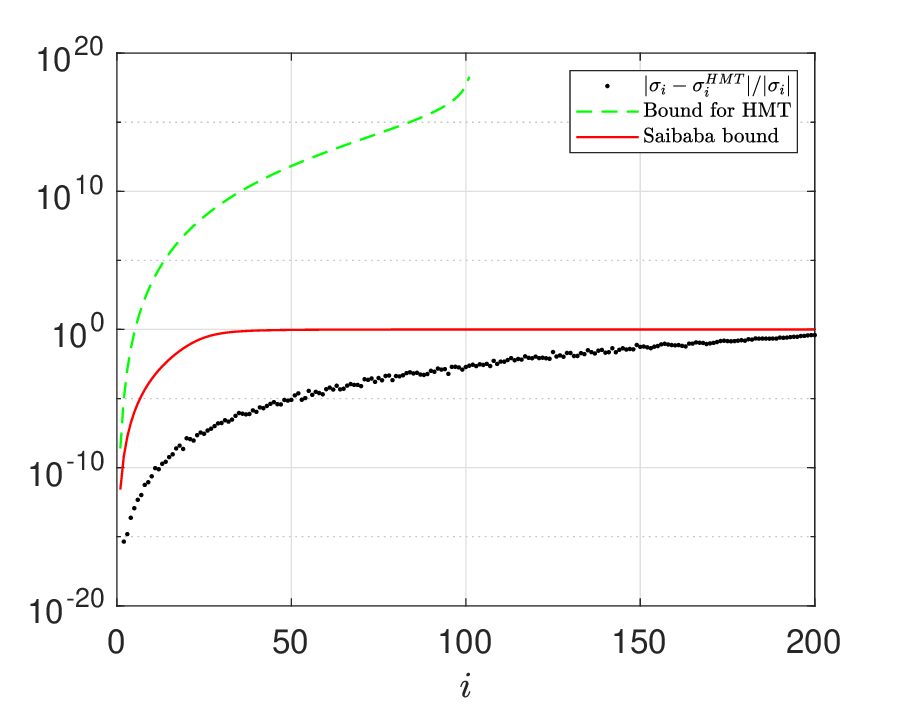}
          \caption{$\tilde{V}$ randomly generated}
    \end{subfigure}
    \caption{\textbf{Comparison between relative bound from \cref{thm:GNbound} and bound in \cite{saibaba2019randomized}.} We show the relative error in extracting the first $r = 200$ singular values by the HMT method with $\ell = 0$ (black dots), Saibaba's bound in \cite{saibaba2019randomized} (red) and the relative bound derived by \cref{thm:GNbound} (dashed green). We considered $A$ with algebraically decaying singular values, and approximate subspaces generated (a) starting from a Gaussian matrix subject to one power iteration and (b) randomly.}
    \label{fig:Saibaba}
    \end{figure}

Figure \ref{fig:Saibaba} shows a comparison of the bound in \cite{saibaba2019randomized} and the bound proposed here. We note that for approximate subspace that are randomly generated, Saibaba's bound is very sharp and describes better the error compared to relative bound derived from \cref{thm:GNbound}. However, when the approximate subspace contains information of the matrix $A$, we see that the proposed bound became sharper, while Saibaba's bound does not clearly depict the relative error.
%%%%%%%%%%%%%%%%%%%%%%%%%%%%%%%%%%%%%%%%%%%%%%%%%%%%%%%%%%%%%%%%%%%%%%%% COMPUTABILITY %%%%%%%%%%%%%%%%%%%%%%%%%%%%%%%%%%%%%%%%%%%%%%%%%%%%%%%%%%%%%%%%%%%%%%%%
\subsection{Computability} \label{subsec:Computability}
 
    A bound is usually derived to be used as a theoretical tool to study the behavior of the bounded quantity and to be used in practical computations. In the previous sections, we have shown how the derived bounds can serve the first purpose. However, as expressed so far, these bounds are impractical to use. Indeed, the presence of the exact singular values makes the bounds useless in practical computations. In this section, we discuss an approach that can overcome this issue. Once again we analyze the \GNa, and, for the sake of simplicity, we restrict ourselves to the case where no oversampling is employed ($\ell = 0$). As seen before, the theory can be easily generalized to the other strategies and to the case with oversampling, see Supplementary Material for the implementation of each case. 
    
    To obtain bound \cref{eq:GNbound}, we used \cref{eq:GNperturbation} to interpret the \GNa{} as a perturbation of $\bA$, i.e. $\bA_{\mtiny{GN}} = \bA - E_{\mtiny{GN}}$. This results in the final expression of the bound including quantities that depend on $\bA$, and consequently on $A$, which are typically unknown in practice. To address this, another less intuitive interpretation is trivially available. Specifically, the roles can be reversed, viewing the matrix $\bA$ as a perturbation of the corresponding \GNa, i.e. $\bA = \bA_{\mtiny{GN}} + E_{\mtiny{GN}}$, where the perturbation $E_{\mtiny{GN}}$ remains the same. We can then apply \Cref{thm:2by2} to this interpretation.
    
\begin{figure}
    \begin{subfigure}{.5\textwidth}
          \centering
          \includegraphics[width=\linewidth]{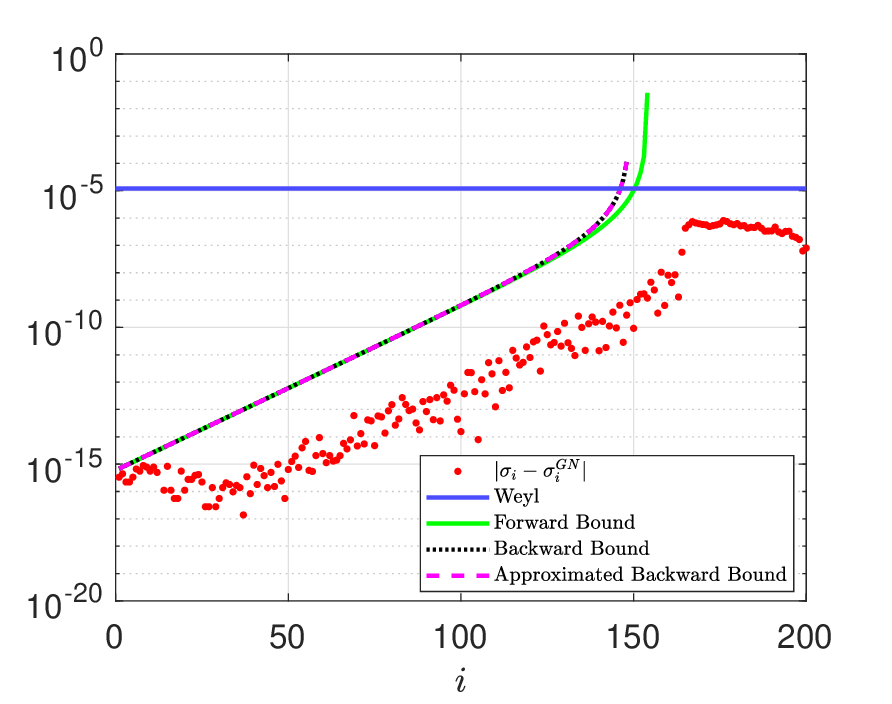}
          \caption{Exponential Decay}
          \label{fig:Computability_exp}
    \end{subfigure}%
    \begin{subfigure}{.5\textwidth}
          \centering
          \includegraphics[width=\linewidth]{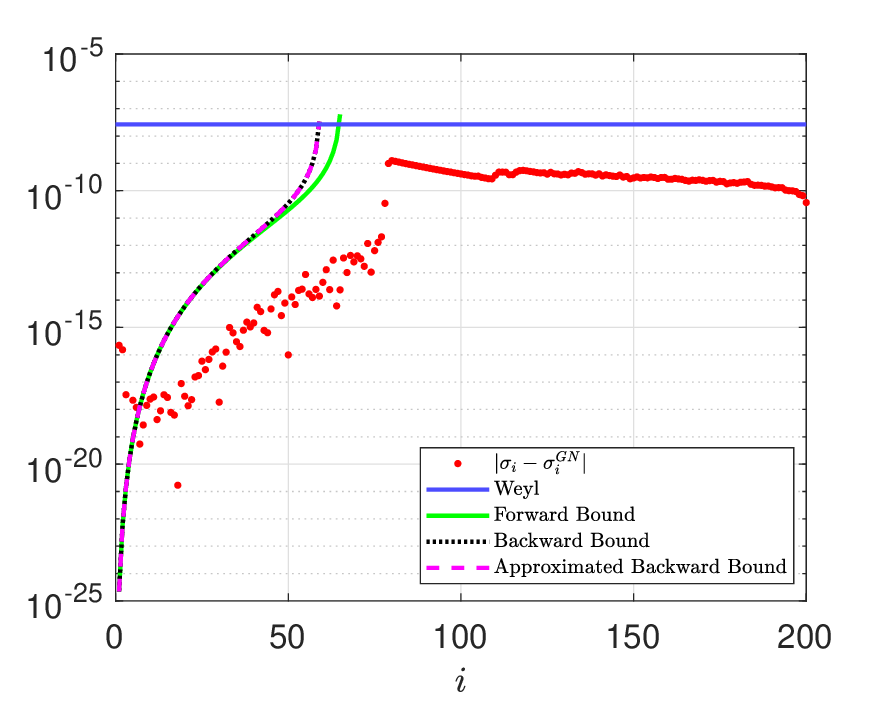}
          \caption{Algebraic Decay}
          \label{fig:Computability_alg}
    \end{subfigure}
    \caption{\textbf{Computable Bound for \GN{} without oversampling.} We show the error in extracting the first $r = 200$ singular values by \GNa{} with $\ell = 0$ (red dots), Weyl’s inequality \cref{eq:Weyl} (blue), the forward bound \cref{eq:GNbound} (green), the backward bound \cref{eq:GNBackwardBound} (black) and the approximated backward bound (magenta). Both exponential (a) and algebraic (b) singular values decays are considered.}
    \label{fig:Computability}
\end{figure}

\begin{corollary}
\label{thm:GNBackwardbound}
    Consider an $m\times n$ matrix $A$ with singular values $\{\sigma_i\}_{i=1}^n$ and approximations $\tV \in \mathbb{R}^{n\times r}$ and $\tU\in \mathbb{R}^{m\times r}$ to its leading singular subspaces. Consider the \GNa\ 
    $A_{\mtiny{ GN,\tV,\tU}} = A\tV (\tU^*A\tV)^\dagger \tU^* A$ and its singular values $\{\sigma^{\mtiny{GN}}_i\}_{i=1}^n$. Define $\bA = Q_1^*AQ_2$, with $Q_{1}= \begin{bmatrix} \tU & \tU_{\perp}\end{bmatrix}$ and $Q_{2}= \begin{bmatrix} \tV & \tV_{\perp}\end{bmatrix}$, and its blocks as in \eqref{def:Abar}. Let, for $i=1, \dots, n$,
    \begin{displaymath}
    \btau_{i} = \frac{ \max\{\|\bA_{12}\|_{2}, \|\bA_{21}\|_2\}}{\displaystyle \min_\minindex |\sigma_{i}(A_{\mtiny{GN}}) - \sigma_\minindex (\bA_{21}\bA_{11}^{\dagger}\bA_{12})| - 2\left\Vert E_{\mtiny{GN}} \right\Vert_{2}}.
\end{displaymath}
Then, for all $i$ for which $\btau_i>0$, it holds
\begin{equation}
\label{eq:GNBackwardBound}
    | \sigma_i - \sigma^{\mtiny{GN}}_i | \leq  \left\Vert  \bA_{22} - \bA_{21}\bA_{11}^\dag \bA_{12} \right\Vert_2 \btau_i^2.
\end{equation}
\end{corollary}
We refer to bound \cref{eq:GNBackwardBound} as the backward bound. Notice that instead of the exact singular values, as in \cref{eq:GNbound}, we now have the computed singular values. However, the gap in the denominator of $\btau_i$ is more complicated to compute. For this reason, we can use the approximation
\begin{displaymath}
\min_\minindex |\sigma_{i}(A_{\mtiny{GN}}) - \sigma_\minindex (\bA_{21}\bA_{11}^{\dagger}\bA_{12})| \approx \sigma_i(A_{\mtiny{GN}}) - \|\bA_{21}\bA_{11}^{\dagger}\bA_{12}\|_2,
\end{displaymath}
to obtain an approximated form of the backward bound. \Cref{fig:Computability} shows that, for both exponentially and algebraically decaying exact singular values, both the backward bound and the approximated backward bound are not far from the forward bound \cref{eq:GNbound}. In practice, the backward bound can be further approximated by e.g. estimating norms by random sampling \cite{martinssonRandomizedNumericalLinear2020a}, approximating the norm of the Schur complement \cite{benziNumericalSolutionSaddle2005a},\cite[Sec. 8.2]{elmanFiniteElementsFast2005}, and, for the \GNa{} with oversampling, making use of properties of the orthogonal projection $\bA_{11}^\dagger \bA_{11}$.
 
%%%%%%%%%%%%%%%%%%%%%%%%%%%%%%%%%%%%%%%%%%%%%%%%%%%%%%%%%%%%%%%%%%%%%%%%%%%%%%%%%%%%%%%%%%%%%%%%%%%%%%%%%%%%%%%%%%%%%%%%%%%%%%%%%%%%%%%%%%%%%%%%
%%%%%%%%%%%%%%%%%%%%%%%%%%%%%%%%%%%%%%%%%%%%%%%%%%%%%%%%%%%%%%%%%%%%%%%% CONCLUSION %%%%%%%%%%%%%%%%%%%%%%%%%%%%%%%%%%%%%%%%%%%%%%%%%%%%%%%%%%%%
%%%%%%%%%%%%%%%%%%%%%%%%%%%%%%%%%%%%%%%%%%%%%%%%%%%%%%%%%%%%%%%%%%%%%%%%%%%%%%%%%%%%%%%%%%%%%%%%%%%%%%%%%%%%%%%%%%%%%%%%%%%%%%%%%%%%%%%%%%%%%%%%
\section{Discussion} \label{sec:Conclusion}
We conclude with a discussion of the possible extensions and open problems.

We believe there is room to further explore how \Cref{thm:2by2} can be used as a tool to formally compare strategies to extract singular values. For example, a study of the dominant term in each bound could lead to further understanding of a possible accuracy hierarchy.

As presented, the technique to improve the bound for the \GNa{} when $\ell>0$ is based on heuristics. Thus, there is space to investigate the formal validity of such strategy and to devise alternative ones.

It is worth noting that the behavior of the presented bounds differs from what is illustrated in \Cref{sec:NumIll} when, instead of good approximations of the leading singular subspaces, we take $\tU$ and $\tV$ to be random matrices. Our experiments show that, in this case, the change of the errors in extracting the leading singular values of the \GNa{} depends less on $i$, giving very similar accuracy errors for the first $r$ singular values. In fact, in this case Weyl's inequality is more descriptive than in the previous experiments, indicating that the need to find better bounds is diminished.

Finally, the presented analysis could be used to characterize the difference in the behavior of the \GNa{} with and without oversampling. This would have an impact also on the analysis of the \GN{} low-rank approximation.

\section*{Acknowledgments}
This work is supported by the the EPSRC Grants No. EP/Y008200/1,  EP/Y010086/1 and 
EP/Y030990/1, and 
the Ada Lovelace Centre Programme at the Scientific Computing Department, STFC.  
\bibliography{main_r}
\bibliographystyle{siamplain}
\end{document}